\documentclass[11pt]{article}
\date{March 31, 2016}

\usepackage {amsmath,amsthm,amsfonts,amssymb,multicol}
\usepackage{graphicx}
\usepackage{booktabs}
\usepackage{array}
\usepackage{subfigure}
\usepackage[font=scriptsize]{caption}
\usepackage{enumerate}
\usepackage[nodayofweek]{datetime}
\usepackage{url}
\usepackage{hyperref}
\usepackage{geometry}
\usepackage[english]{babel}
\usepackage{url}

\usepackage{verbatim} 
\usepackage{amsthm} 
\usepackage{enumitem}
\usepackage{enumerate}
\usepackage{bbm} 
\usepackage[normalem]{ulem} 
\usepackage{bm}
\usepackage{fullpage}
\usepackage{authblk}

\newcommand{\be}{\begin}
\newcommand{\e}{\end}
\newcommand{\beq}{\begin{equation}}
\newcommand{\eeq}{\end{equation}}
\newcommand{\beqs}{\begin{equation*}}
\newcommand{\eeqs}{\end{equation*}}
\newcommand{\bal}{\begin{align}}
\newcommand{\eal}{\end{align}}
\newcommand{\bals}{\begin{align*}}
\newcommand{\eals}{\end{align*}}

\newcommand{\ol}{\overline}

\newcommand{\R}{\mathbb{R}}
\renewcommand{\S}{\mathbb{S}}
\newcommand{\C}{\mathbb{C}}

\renewcommand{\l}{\left}
\renewcommand{\r}{\right}

\renewcommand{\d}{\mathrm{d}} 

\newcommand{\set}[1]{\mathbb{#1}}

\newcommand{\curly}[1]{\mathcal{#1}}

\newcommand{\setof}[2]{\left\{ #1\; : \;#2 \right\}}
\newcommand{\cc}{\Subset}

\newcommand{\om}{\omega}
\newcommand{\Om}{\Omega}
\newcommand{\eps}{\varepsilon}
\newcommand{\lam}{\lambda}

\newcommand{\gam}{\gamma}
\newcommand{\Gam}{\Gamma}

\newcommand{\al}{\alpha}
\newcommand{\de}{\delta}
\newcommand{\De}{\Delta}
\newcommand{\Del}{\Delta}

\newcommand{\vp}{\varphi}

\newcommand{\scp}[2]{\langle#1,#2\rangle}

\newcommand{\jap}[1]{\l\langle #1\r\rangle}

\renewcommand{\it}{\infty}
\newcommand{\IT}{\infty}

\newcommand{\supp}{\,\textnormal{supp}}

\newcommand{\del}{\partial}

\newcommand{\dx}{\,\d x}
\newcommand{\dy}{\,\d y}

\newcommand{\ft}[1]{\widehat{{#1}}}			

\theoremstyle{definition}

\theoremstyle{remark}

\def\dotuline{\bgroup
  \ifdim\ULdepth=\maxdimen  
   \settodepth\ULdepth{(j}\advance\ULdepth.4pt\fi
  \markoverwith{\begingroup
  \advance\ULdepth0.08ex
  \lower\ULdepth\hbox{\kern.15em .\kern.1em}%
  \endgroup}\ULon}

\def\dashuline{\bgroup
  \ifdim\ULdepth=\maxdimen  
   \settodepth\ULdepth{(j}\advance\ULdepth.4pt\fi
  \markoverwith{\kern.15em
  \vtop{\kern\ULdepth \hrule width .3em}%
  \kern.15em}\ULon}

\begin{document}

\title{On the H\"older regularity for the fractional Schr\"o\-dinger equation and its improvement for radial data}

\author{Marius Lemm\thanks{mlemm@caltech.edu} }
\affil{Mathematics Department, Caltech}

\maketitle

\abstract{
We consider the linear, time-independent fractional Schr\"odinger equation
$$
    (-\Del)^s \psi+V\psi=f\quad \text{on } \Om\subset \set{R}^N.
$$
We are interested in the local H\"older exponents of distributional solutions $\psi$, assuming local $L^p$ integrability of the functions $V$ and $f$. By standard arguments, we obtain the formula $2s-N/p$ for the local H\"older exponent of $\psi$ where we take some extra care regarding endpoint cases. For our main result, we assume that $V$ and $f$ (but not necessarily $\psi$) are radial functions, a situation which is commonplace in applications. We find that the regularity theory ``becomes one-dimensional'' in the sense that the H\"older exponent improves from $2s-N/p$ to $2s-1/p$ away from the origin. Similar results hold for $\nabla\psi$ as well.
}

\section{Introduction}

The fractional Laplacian $(-\De)^s$ appears in the modeling of many real world phenomena with anomalous diffusion, see the survey paper \cite{Hitchhiker} for an extensive list of applications. In physics, the closely related operator $\sqrt{-\De+m^2}-m$ ($m>0$) is used to describe the kinetic energy of particles in a pseudorelativistic regime.

From an analytical perspective, $(-\De)^s$ is a prototypical \emph{non local} operator. It is then a natural challenge to establish local regularity results for solutions $\psi$ to equations involving $(-\De)^s\psi$. We call a regularity result ``local'' if the assumption that the data (e.g.\ the functions $V$ and $f$ in \eqref{eq:FSE0}) is nice (say in a $C^k$ or $L^p$ sense) on some ball $B\subset \R^N$, leads to certain regularity of the solution $\psi$ inside of $B$.

Arguably the simplest equation which involves the fractional Laplacian is the ``fractional Poisson equation'' $(-\De)^s\psi=f$. On the whole space $\R^N$, this equation is solved by the Riesz potential of $f$. In this paper, we study the linear fractional Schr\"odinger equation 
\beq
\label{eq:FSE0}
    (-\Del)^s \psi+V\psi=f.
\eeq
which can be understood as a lower order perturbation of the fractional Poisson equation from the perspective of pseudo differential operators (the principal symbol $|\xi|^{2s}$ is unchanged). The addition of a somewhat singular potential $V$ (e.g.\ the Coulomb potential) is the kind of perturbation that is most relevant to applications in quantum mechanics \cite{CMS, Herbst, LiebYau}. The main point of this paper is to study the local regularity of solutions $\psi$ to the concrete equation \eqref{eq:FSE0} \emph{quantitatively}. That is, we are interested in obtaining the best possible value of the local H\"older exponent of $\psi$ and we place particular focus on the case when $V$ and $f$ are radial.
 
 We note that our perspective is different from one that is commonly taken in the PDE community, where regularity is proved in the more general case when the highest order term $(-\De)^s\psi$ is perturbed but with less focus on quantitative information (e.g.\ on the best possible value of the H\"older exponent). For second-order differential equations, this is an old topic with the key contributions coming from de Giorgi, Nash and Moser \cite{deGiorgi,Moser, Nash} for equations in divergence form and Krylov and Safonov \cite{KS1, KS2} for equations in non divergence form. For non local operators, the topic has blossomed in recent years, see e.g.\ \cite{BassKassmann05, BassLevin02, BCF12, CaffarelliStinga14,  CS07, CS09, CS14,FelmerQuaasTan12, JS15, KRS14, Ros-OtonSerra14, Silvestre07}. 

We were motivated in part by the papers \cite{DallAcquaFournaisSorensenStockmeyer, DallAcquaFournaisSorensenStockmeyer0} where analyticity of atomic eigenfunctions away from the nuclei is shown in the pseudorelativistic regime. In terms of methods, we draw on a helpful localization trick due to Silvestre \cite{Silvestre07} which was further developed by Felmer, Quaas and Tan \cite{FelmerQuaasTan12}.

\paragraph{Main results.}
We now describe the main results in words. The precise statements are in Section 3. Overall, the present paper is devoted to local and quantitative H\"older regularity results for solutions $\psi$ to \eqref{eq:FSE0}.

 Our \textbf{first result}, Theorem \ref{thm:Holder1}, concerns \eqref{eq:FSE0} on the whole space $\R^N$. It says that if $V,f\in L^p(B)$ on some ball $B\subset \R^N$, then any solution $\psi$ is H\"older continuous inside of $B$ and its H\"older exponent is $2s-N/p$ whenever this number lies in $(0,1)$. Moreover, if $2s-N/p$ equals one, then $\psi$ is Lipschitz continuous up to a logarithmic correction factor. A similar result holds for $\nabla\psi$ with $2s-N/p$ replaced by $2s-N/p-1$. The result is obtained by the standard bootstrap argument and a localization trick due to \cite{Silvestre07}. Our contribution here is mainly to take some care to include the endpoint cases (a) $2s=N$ (where the inverse of $(-\De)^s$ is given by the logarithmic potential, also for $s=1/2$ and $N=1$) and (b) $2s-N/p=1$ where we prove Lipschitz continuity up to a logarithmic correction factor. This part has some relation to the study of global mapping properties of Riesz kernels on $\R^N$, a classical topic which goes back at least to the early 60s \cite{Landkof, SteinWeiss, Yudovic61} and is still being studied \cite{GargSpector14, GargSpector15}.

Our \textbf{main contribution}, Theorem \ref{thm:radial}, is the following: Assume that the data $V$ and $f$ are \emph{radial} functions, a situation which is commonplace in applications. Then, in the formulae $2s-N/p$ (or $2s-N/p-1$) for the H\"older exponents, $N$ can be replaced by \emph{one} away from the origin. In other words, the \emph{effective dimension governing regularity away from the origin is equal to one due to the radial symmetry}. We emphasize that $\psi$ is not assumed to be radial here. The idea of an effective dimensional reduction (of the regularity theory) in the presence of symmetry is intuitive but it seems that it has not been emphasized so far in the context of non local operators. (For the ordinary Laplacian it is essentially trivial, as we discuss below.) 


We also prove two \textbf{further results}.
\be{itemize}
\item The results of Theorem \ref{thm:Holder1} can be extended to the case when the equation holds on a finite-measure domain $\Om\subset \R^N$, see Theorem \ref{thm:Holder1local}. This requires a different localization argument, a non local variant of the Leibniz rule, and the observation that the corresponding localization error can never worsen the local integrability of $\psi$ (Lemma \ref{lm:locerror}). The non local Leibniz rule is related to the ``non local IMS formula'' which first appeared in \cite{LiebYau} and is ascribed to Michael Loss there. 
 \item The H\"older exponents derived in Theorems \ref{thm:Holder1} and \ref{thm:Holder1local} are optimal, see Proposition \ref{thm:optimality}.
\e{itemize}

Throughout the paper, we work with a new and relatively weak notion of distributional solution, see Definition \ref{defn:distributional}. In particular, it is weaker than the usual assumption that a solution lies in the Sobolev space $H^s(\R^N)$.

\paragraph{Radial symmetry for $s=1$.}
To understand the heuristics of the dimensional reduction in Theorem \ref{thm:radial}, it is instructive to consider the case of the ordinary Laplacian, where it is trivial. Let $N\geq 1$ and let $V:\R^N\to \C$ be radial. Recall that $L^2(\R^N)$ may be decomposed into the subspaces $\curly{H}_l$ spanned by the spherical harmonics of order $l$ (defined in \eqref{eq:Hldefn}). Each $\curly{H}_l$ is invariant under the action of the Laplacian and under multiplication by $V$. Hence, we may restrict to solutions $\psi\in\curly{H}_l$ which solves $-\De\psi+V\psi=0$ ($f$ only appears when $l=0$ and we ignore this possibility here). In spherical coordinates $x=r\om$ with $r>0$ and $\om\in \S^{N-1}$, we can decompose any function $\psi\in\curly{H}_l$ as $\psi(x)=\Psi(r) Y_l(\om)$ where $Y_l$ is smooth. To find the equation satisfied by $\Psi$ we write the Laplacian in spherical coordinates
 \beq
 \label{eq:radialLaplace}
   -\De \psi=-\frac{1}{r^{N-1}}\del_r \l(r^{N-1}\del_r \psi\r)  -\frac{1}{r^2} \De_{\set{S}^{N-1}}\psi.
 \eeq
 Here $\De_{\set{S}^{N-1}}$ is the Laplace-Beltrami operator on the sphere and we have $\De_{\set{S}^{N-1}} Y_l=\lam_lY_l$ for some $\lam_l\in\R$. This yields the second-order ODE satisfied by $\Psi$ on $\R_+$. Note that, away from the origin, all powers of $r$ are smooth functions. In particular, the second term in \eqref{eq:radialLaplace} amounts to a smooth, bounded potential, which can just be added to the more singular potential $V$ without changing its local $L^p$ properties. Moreover, if we change the unknown from $\Psi(r)$ to $u(r):=r^{(N-1)/2}\Psi(r)$, a short computation shows that $u$ satisfies a \emph{one-dimensional Schr\"odinger equation} with a potential that still has the same local $L^p$ regularity as $V$. Hence, the one-dimensional regularity theory \cite{GilbargTrudinger, LiebLoss} applies to $u$. Since $r^{(N-1)/2}$ and $Y_l$ are smooth functions away from the origin, $\psi$ has the \emph{same} regularity as $u$ and so the effective dimension governing the regularity of $\psi$ has indeed been reduced to one. 

This argument, which was based on \eqref{eq:radialLaplace}, breaks down completely for the fractional Laplacian. The main point of this paper is to prove that the conclusion holds nonetheless.

\subsection*{Organization of the paper}

In \textbf{Section 2}, we give the relevant definitions, in particular we provide a new notion of distributional solution for the non local equation \eqref{eq:FSE0}. Then we characterize all distributional solutions to the fractional Poisson equation $(-\De)^s\psi=f$ on $\set{R}^N$ in terms of the Riesz potential, see Proposition \ref{prop:riesz}. 

In \textbf{Section 3}, we state our \emph{main results}. The standard regularity theory on $\R^N$ is discussed in Theorem \ref{thm:Holder1}. Our main contribution is Theorem \ref{thm:radial}, on the improvement of the H\"older exponents in the presence of spherical symmetry. Then we give two further results: (a) Theorem \ref{thm:Holder1local} which generalizes Theorem \ref{thm:Holder1} to finite-measure domains $\Om\subset\R^N$ and (b) Proposition \ref{thm:optimality} which notes the optimality of the derived H\"older exponents. 

In \textbf{Section 4}, we prove Theorem \ref{thm:Holder1} by the standard bootstrap argument, combined with Silvestre's localization trick 

In \textbf{Section 5}, we prove Theorem \ref{thm:radial} using ideas from the previous section, properties of hypergeometric functions and an amusing trick due to Frank and Lenzmann which ``trades off angular momentum for dimension'' (in physics jargon). 

In \textbf{Section 6}, we prove Theorem \ref{thm:Holder1local}. The idea is to extend the Schr\"odinger equation from $\Om$ to $\R^N$ and to apply Theorem \ref{thm:Holder1}. Technically, this requires bounds on the ``localization error'' coming from a non local version of the Leibniz rule, see Lemma \ref{lm:locerror}, which may be of independent interest. 


\section{Preliminaries}
 \subsection{Definitions}
\paragraph{The fractional Laplacian and the Riesz kernel. }
 Let $0<s\leq 1$. For any Schwartz function $\vp$, the \emph{fractional Laplacian} is defined as the Fourier multiplier
\beqs
    (-\De)^s \vp = (|\cdot|^{2s} \ft{\vp})^{\vee}
\eeqs
where we use the following convention for the Fourier transform and its inverse
\beqs
    \ft{\vp}(\xi)=(2\pi)^{-N/2} \int_{\set{R}^N} e^{-i \xi x} \vp(x)   \d x,\qquad 
    \check \vp(x)= (2\pi)^{-N/2} \int_{\set{R}^N} e^{i  \xi x} \vp(\xi)   \d \xi.
\eeqs
When $s<1$, it is well-known that one can use Plancherel's identity to rewrite $(-\De)^s$ as the singular integral operator
\beq
\label{eq:FLdefn}
		(-\Del)^s \vp(x) = C_{N,s}\,  \mathrm{P.V.}\int \frac{\vp(x)-\vp(y)}{|x-y|^{N+2s}}\dy,\
\eeq
where
\beqs
    C_{N,s}= \pi^{-N/2} 2^{2s} \frac{\Gam\l(\frac{N+2s}{2}\r))}{|\Gam(-s)|}
\eeqs
Here, $\mathrm{P.V.}$ stands for the principal value. Note that it can be dropped when $s<1/2$. We define the (slightly generalized) \emph{Riesz kernels} $k_{2s}$ on $\R^N$ by
\begin{align}
\label{eq:rieszdefn}
	k_{2s}(x):=
    \be{cases}
        D_{N,s} |x|^{-(N-2s)},\quad &\text{if } 2s\neq N\\
        \om_{N-1}^{-1} \log(|x|^{-1}), \quad &\text{if } 2s=N=2\\
        -\sqrt{\frac{2}{\pi}}(\gam_{EM}+\log|x|)\quad &\text{if } 2s=N=1
    \e{cases}
\end{align}
Here $\om_{N-1}$ is the measure of the $(N-1)$-sphere, $\gam_{EM}\approx 0.577$ is the Euler-Mascheroni constant and
 \beqs
 D_{N,s}= \frac{1}{\pi^{N/2} 2^{2s}} \frac{\Gam\l(\frac{N-2s}{2}\r)}{\Gam(s)}.
\eeqs

\paragraph{Function spaces. }
Our results will say that solutions $\psi$ lie in the following subspaces  of continuous functions. They are just the usual H\"older spaces when the H\"older exponent $\beta<1$, but they acquire a logarithmic correction factor to Lipschitz continuity when $\beta=1$. 

\be{defn}[H\"older spaces]
\label{defn:holder}
Let $U\subset\set{R}^N$ be open, let $g:U\to \set{C}$ and $0<\beta\leq 1$. We say
\beqs
    g\in \tilde C^{0,\beta}(U)
\eeqs
if there exists a constant $c>0$ such that for all $x,y\in U$, it holds that
\beqs
    |g(x)-g(y)|\leq c
    \be{cases}
     |x-y|^\beta,\quad &\text{if } \beta<1\\
     |x-y| \log\l(|x-y|^{-1}\r), \quad &\text{if } \beta=1.
    \e{cases}
\eeqs
Moreover, we say $g\in \tilde C^{1,\beta}(U)$, if $g\in C^1(U)$ and $\nabla g\in \tilde C^{0,\beta}(U)$.
\e{defn}

Throughout, we work with complex-valued functions, since we have applications to quantum mechanics in mind.
To state our results under rather sharp integrability assumptions, we make

\be{defn}
\label{defn:weightedL1}
Let $\jap{x}:=(1+x^2)^{1/2}$. For $\beta\in \R$, we define the weighted spaces
\beqs
    \jap{\cdot}^{\beta}L^1(\set{R}^N) := \setof{f:\set{R}^N\rightarrow\set{C}}{\int_{\set{R}^N} \jap{x}^{-\beta} |f(x)|\d x<\it}.
\eeqs
Moreover, for $\gam \geq 1$, we define the following combination of Lebesgue spaces
$$
\curly{L}_\gam(\R^N):=
\be{cases}
\bigcup\limits_{0<\de<\gam-1} \l(L^{1+\de}(\R^N)+L^{\gam-\de}(\R^N)\r),&\text{ if }\gam>1,\\
\bigcup\limits_{\de>0}\l(L^{1+\de}(\R^N)\cap \jap{\cdot}^{-\de}L^1(\R^N)\r), &\text{ if }\gam=1.
\e{cases}
$$
\e{defn}


A sizable subset of functions in $\curly{L}_{\gam}(\R^N)$ are those for which there exists $\de>0$ such that they lie in $L^{1+\de}_{\mathrm{loc}}(\R^N)$ and decay like $\jap{x}^{-N/\gam-\de}$ at infinity. We will usually take $\gam=N/(2s)$. The definition is designed such that $f\in\curly{L}_{N/(2s)}(\R^N)$ guarantees that
$$\int_{\R^N} ((-\De)^s\vp)(x) (k_{2s}*f)(x)\d x<\it$$ holds for any $\vp\in C_0^\it(\R^N)$ by the Hardy-Littlewood-Sobolev inequality, see e.g.\ the proof of Proposition \ref{prop:riesz} in Appendix \ref{app:propriesz}.

\paragraph{Conventions.} 
We write $C,C',\ldots$ for positive constants which depend only on the values of fixed parameters, e.g.\ the dimension $N$. The numerical value of $C,C',\ldots$ may change from line to line. We write $p'$ for the H\"older dual of $p\in [1,\it]$, i.e.\ $p'=p/(p-1)$ with the convention that $1/\it=0$. We denote duality products by $\scp{\cdot}{\cdot}$. For two open sets $A,B\subset\set{R}^N$, we write $
 A\cc B$ to express that $\ol{A}\subset B$. We use the term ``cutoff function'' for a function in $C_0^\it(\R^N)$ which takes values only in $[0,1]$. An ``annulus'' in $\R^N$ is always assumed to be centered at the origin.

\emph{Throughout the entire paper, $N\geq 1$ and $0<s\leq\min\{1,N/2\}$.} 

\paragraph{The case $1<s<N/2$.}
The techniques in this paper generalize in principle to all $0<s<N/2$. This is because the Riesz kernel $k_{2s}$ is locally integrable for all $0<s<N/2$. However, there is not a nice singular integral representation like \eqref{eq:FLdefn} for such $s$ and alternatives (see e.g.\ p.\ 46f in \cite{Landkof}) would lead to some case distinctions for the integrability conditions on distributional solutions. Therefore we restrict to $0<s<1$ for the sake of simplicity of presentation. Regularity questions for the full range $0<s<N/2$ are addressed e.g.\ in the related papers \cite{JLX, L}.

\subsection{Distributional solutions}
Since $(-\Del)^s$ is a non local operator when $0<s<1$, we need to modify the usual definition of a distributional solution.
\be{defn}[Distributional solutions]
\label{defn:distributional}
Let $\Om\subset\set{R}^N$ be open. Let $V,f:\set{R}^N\rightarrow \set{C}$ with $f\in L^1_{\mathrm{loc}}(\Om)$. We say that $\psi$ is a \emph{distributional solution} of the fractional Schr\"odinger equation
\beq
\label{eq:FSE}
		(-\Del)^s \psi + V\psi=f \quad \text{on } \Om,
\eeq
if it satisfies
\be{itemize}
\item $\psi\in \jap{\cdot}^{N+2s} L^1(\set{R}^N)$, when $s<1$, or $\psi\in  L^1_{\mathrm{loc}}(\Om)$, when $s=1$,
\item $V\psi\in L^1_{\mathrm{loc}}(\Om)$,
\item for all $\vp\in C_0^\it(\Om)$, it holds that
\beq
\label{eq:FSEdistribution}
		\scp{(-\Del)^s\vp}{\psi} + \scp{\vp}{V\psi} = \scp{\vp}{f}.
\eeq
\e{itemize}
\e{defn}

\be{prop}
\label{prop:finite}
Under the integrability assumptions in Definition \ref{defn:distributional}, the expression in \eqref{eq:FSEdistribution} makes sense.
\e{prop}

This is proved in Appendix \ref{app:propfinite}. We note that the assumption $\psi\in \jap{\cdot}^{N+2s} L^1(\R^N)$ differs from the usual assumption $\psi\in L^1_{\mathrm{loc}}(\Om)$ for distributional solutions to differential equations. It is necessary to control the term $\scp{(-\Del)^s\vp}{\psi}$ in \eqref{eq:FSEdistribution} at infinity. (Note that $(-\Delta)^s \vp$ may be non-zero outside of $\supp\,\vp$ and it decays with a certain polynomial rate, see \eqref{eq:FLdefn}.)

\be{rmk}
We have $H^s(\R^N)\subset \jap{\cdot}^{N+2s} L^1(\set{R}^N)$ by the Sobolev embedding theorem and H\"older's inequality. In particular, our Definition \ref{defn:distributional} of a distributional solution generalizes the common notion of a weak solution in $H^s$.
\e{rmk}

\subsection{Solving the fractional Poisson equation on $\set{R}^N$}
An important intermediate result for us is to explicitly solve the fractional Poisson equation (i.e.\ \eqref{eq:FSE} with $V\equiv0$) on $\set{R}^N$ in terms of the Riesz kernel $k_{2s}$ defined in \eqref{eq:rieszdefn} above. Recall our standing assumption that $0<s\leq\min\{1,N/2\}$ and Definition \ref{defn:weightedL1} of the function space $\curly{L}_\gam$.

\be{prop}
\label{prop:riesz}
Let $f\in \curly{L}_{N/(2s)}(\R^N)$ and suppose that $\psi$ solves
\beqs
    (-\De)^s\psi = f \quad \text{on } \set{R}^N
\eeqs
in the sense of Definition \ref{defn:distributional}. Then:
\be{enumerate}[label=(\roman*)]
\item If $0<s\leq 1/2$, there exists $w\in\set{C}$ such that, for a.e.\ $x\in\set{R}^N$,
\beqs
    \psi(x)= w+(k_{2s} * f)(x).
\eeqs

\item If $1/2<s\leq 1$, there exist $w,w_1,\ldots,w_N\in\set{C}$ such that, for a.e.\ $x\in\set{R}^N$,
\beqs
    \psi(x) = w+ \sum_{j=1}^N w_j x_j+(k_{2s}*f)(x).
\eeqs
\e{enumerate}
\e{prop}

\be{rmk}
\be{enumerate}[label=(\roman*)]
\item The fact that the Riesz potential \eqref{eq:rieszdefn} solves the fractional Poisson equation is well-known, see e.g.\ Landkof's book \cite{Landkof}. Some extra work is required here because we identify all solutions and we work with a rather weak notion of a distributional solution, which may even grow at infinity. 
\item We see that distributional solutions to the fractional Poisson equation are only unique up to affine functions. Of course, these do not affect regularity.
\item If $s=1$, it suffices to assume that $\psi\in L_{\mathrm{loc}}^1(\R^N)$. This is also true for the main results stated in the next section, but we will not distinguish the $s=1$ case hereafter to keep the presentation simple.
\e{enumerate}
\e{rmk}

The proof of Proposition \ref{prop:riesz} for $2s<N$ is based on the following well known fact, see e.g.\ Lemma 2 on p.\ 117 in \cite{SteinSIDPF}. It may be intuitively understood from the rotational invariance and identical scaling behavior of both functions.

\be{prop}
\label{prop:rieszkernel}Let $2s<N$. Then $\ft{k_{2s}}(\xi) = (2\pi)^{-N/2}|\xi|^{-2s}$, where the Fourier transform is taken in the sense of tempered distributions.
\e{prop}

The proof of Proposition \ref{prop:riesz} consists of approximation arguments that justify taking the Fourier transform to apply this fact. Tempered distributions (the dual space of the Schwartz functions) are the natural subspace of the distributions to define the Fourier transform on, but some complications arise because $(-\De)^s$ does not map Schwartz functions to Schwartz functions. The details are in Appendix \ref{app:propriesz}.

\section{Main results}
We recall our standing assumption that $0<s\leq\min\{1,N/2\}$ and we recall our convention that $1/\it=0$ and that $A\cc B$ means $\ol{A}\subset B$. The function spaces $\tilde C^{0,\al}$ and $\curly{L}_\gam$ were introduced in Definitions \ref{defn:holder} and \ref{defn:weightedL1}.

\subsection{H\"older and almost Lipschitz continuity on $\R^N$}
As discussed before, our first main result establishes that there is a \emph{local} regularity theory for the linear fractional Schr\"odinger equation on $\set{R}^N$, in the sense that local $L^p$ properties of $V$ and $f$ lead to local H\"older continuity of the solution $\psi$. 

\be{thm}
\label{thm:Holder1}
Let $B\subset \set{R}^N$ be an open ball such that $f,V\in L^p(B)$ for some $ 1< p \leq \it$. Let $\psi$ be a solution to
 \beqs
		(-\Del)^s \psi + V\psi=f \quad \text{on } \set{R}^N,
\eeqs
 in the sense of Definition \ref{defn:distributional}.
   Suppose that $f,V\psi \in \curly{L}_{N/(2s)}(\R^N)$. Let $B'$ be any open ball such that $B'\cc B$ and set $\al:=2s-N/p$.    Then:

\be{enumerate}[label=(\Roman*)]
\item Let $0<s\leq 1/2$. When $p>\frac{N}{2s}$, it holds that $\psi\in \tilde C^{0,\al}(B')$.

\item Let $1/2<s\leq 1$.

\be{enumerate}
\item[(II.A)] When $\frac{N}{2s}<p\leq \frac{N}{2s-1}$, it holds that $\psi\in \tilde C^{0,\al}(B')$.
\item[(II.B)] When $p>\frac{N}{2s-1}$, it holds that $\psi\in \tilde C^{1,\al-1}(B')$.
\e{enumerate}
\e{enumerate}
\e{thm}

\be{rmk}
\be{enumerate}[label=(\roman*)]
\label{rmk:c}
\item The bounds on $p$ are exactly such that the H\"older exponent lies in the interval $(0,1]$. When the H\"older exponent reaches one (e.g.\ when $2s-N/p=1$), we prove that $\psi$ is Lipschitz continuous up to a logarithmic correction factor. This kind of ``almost Lipschitz continuity'' is common in the endpoint case, see e.g.\ the related works \cite{GargSpector14, MuehlemannThesis}.
\item The key assumption is $f,V\in L^p(B)$ and $p$ enters in the formula for the H\"older exponent $\al$. The additional assumption that $f,V\psi \in \curly{L}_{N/(2s)}(\R^N)$ is rather weak and only slightly stronger than $\psi$ being a distributional solution in the sense of Definition \ref{defn:distributional}. It is present for two reasons: (a) to obtain a good solution theory by ensuring integrability against the Riesz kernel via the Hardy-Littlewood-Sobolev inequality, see Proposition \ref{prop:riesz}, and (b) to control non local tails when using Lemma \ref{lm:step1} (whose proof also uses Proposition \ref{prop:riesz}).
\item For $2s\leq N$, there is no Sobolev embedding theorem into H\"older spaces. Hence, even if one restricts to solutions $\psi \in H^s(\R^N)$, the statement that $\psi\in \tilde C^{0,\al}$ with $\al>0$ is non-trivial.
\e{enumerate}
\e{rmk}

\paragraph{Proof strategy.}
The proof of Theorem \ref{thm:Holder1} is split into the usual two main steps, well known from the regularity theory for the ordinary Laplacian:
\be{itemize}
    \item In Lemma \ref{lm:step1}, we show that $\psi$ is bounded on a ball
      containing $B'$ by a bootstrap argument.
    \item We then use the smoothing properties of the Riesz kernel $k_{2s}$ to conclude the claimed H\"older continuity from the previously shown boundedness of $\psi$.
\e{itemize}
A key ingredient to make these steps work nicely in the non local context, is the localization-on-the-right trick, due to \cite{Silvestre07} and formulated here as Lemma \ref{lm:pwise}. 

\subsection{Improved regularity for radial $V$ and $f$}
We say a function $F:\R^N\to\C$ is ``radial'' iff $F(x)=F_0(|x|)$ for some $F_0:\R_+\to\C$ and we often abuse notation and identify $F\equiv F_0$. 

The following result says that when $V$ and $f$ are radial functions, the H\"older exponent $\al=2s-N/p$ from Theorem \ref{thm:Holder1} can be replaced by $\al_1=2s-1/p$ away from the origin. The intuition for this is that the spherical symmetry reduces the effective dimension governing regularity from $N$ to \emph{one}.

Before we state the result, we briefly review the decomposition of $L^2(\set{R}^N)$ into sectors of spherical harmonics. By Lemma 2.18 in Chapter IV of  \cite{SteinWeiss}, we have the orthogonal decomposition
\beq
\label{eq:L2decomposition}
    L^2(\set{R}^N)=\bigoplus_{l=0}^\it \curly{H}_l(\set{R}^N).
\eeq
Here, using spherical coordinates $x=r\om$ with $r>0$ and $\om\in\set{S}^{N-1}$, the subspace $\curly{H}_l(\set{R}^N)\subset L^2(\set{R}^N)$ is defined by
\beq
\label{eq:Hldefn}
\curly{H}_l:=\mathrm{span}\setof{\Psi(r)Y_l(\om)}{\Psi\in L^2(\R_+; r^{N-1}\d r),\, Y_l \textnormal{ a spherical harmonic of degree $l$}}.
\eeq
For a set $A\subset \R^N$, we define $\min A:=\min\setof{|y|}{y\in A}$. Recall that by convention any ``annulus'' is always assumed to be centered at the origin.

\be{thm}
 \label{thm:radial}
 Let $N\geq 2$. Let $A\subset\set{R}^N$ be an open annulus with $\min A>0$. 
 Let $V,f$ be radial functions satisfying $V,f\in L^p(A)$ for some $1< p\leq \it$. Suppose $\psi\in \curly{H}_l(\set{R}^N)$ solves
 \beq
 \label{eq:FSErad}
 (-\De)^s\psi+V\psi=f,\quad \text{on }\set{R}^N
  \eeq
 in the sense of Definition \ref{defn:distributional} for some $l\geq0$. Suppose that
 \beq
 \label{eq:reg}
 |\cdot|^{-l}V\psi\,\in 
 \be{cases}
 \curly{L}_{1/(2s)}(\R_+; r^{N+2l-1}\d r),\quad &\text{if } 0<s\leq 1/2\\
 \jap{\cdot}^{N+2l-2s}L^1(\R_+; r^{N+2l-1}\d r),\quad &\text{if } 1/2<s\leq 1
 \e{cases}
 \eeq
 and that the same assumption with $l=0$ holds for $f$.
 Let $A'$ be any open annulus such that $A'\cc A$ and set $\al_1=2s-1/p$. Then, in short, the conclusion of Theorem \ref{thm:Holder1} holds with $B'$ replaced by $A'$ and $\al$ replaced by $\al_1$. Explicitly:
\be{enumerate}[label=(\Roman*)]
\item Let $0<s\leq 1/2$. When $p>\frac{1}{2s}$, it holds that $\psi\in \tilde C^{0,\al_1}(A')$.

\item Let $1/2<s\leq 1$.

\be{enumerate}
\item[(II.A)] When $\frac{1}{2s}<p\leq \frac{1}{2s-1}$, it holds that $\psi\in \tilde C^{0,\al_1}(A').$
\item[(II.B)] When $p>\frac{1}{2s-1}$, it holds that $\psi\in \tilde C^{1,\al_1-1}(A').$
\e{enumerate}
\e{enumerate}
\e{thm}

\be{rmk}
\label{rmk:open}
 \be{enumerate}[label=(\roman*)]
 \item The case  $s=1$ is much simpler and a proof was sketched in the introduction.
\item If \eqref{eq:FSErad} holds for $\psi\in\curly{H}_l$ and $l>0$, then necessarily $f\equiv 0$. This follows directly from the orthogonality of the decomposition \eqref{eq:L2decomposition}.
\item The regularity assumption \eqref{eq:reg} replaces the assumption in Theorem \ref{thm:Holder1} that $V\psi,f\in \curly{L}_{N/(2s)}(\R^N)$, see also Remark \ref{rmk:c} (ii). (The function spaces appearing in it are the natural generalizations of the ones from Definition \ref{defn:weightedL1} to more general measure spaces.) Due to the nature of the $l$ dependence, it is essentially a decay assumption on $V\psi$ which becomes stronger as $l$ grows. It is possible to weaken this assumption by multiplying $\psi$ by a cutoff function and using the techniques from Section 6, but we do not dwell on this further.
 \item $\curly{H}_0(\set{R}^N)$ is the subspace of radial functions in $L^2(\set{R}^N)$ and so Theorem \ref{thm:radial} applies to these in particular. In fact, we use an amusing trick to reduce the case of arbitrary $l$ to the radial case $l=0$ by increasing the dimension (from $N$ to $N+2l$, compare \eqref{eq:reg}). This useful (but not strictly necessary) trick is an unpublished observation of R.~Frank and E.~Lenzmann and we are grateful to them.
\item The theorem holds for solutions in $\curly{H}_l(\set{R}^N)$ with $l\geq 0$ arbitrary. Of course it immediately generalizes to all finite combinations of such $\psi$. However, the bounds do not decay sufficiently fast in $l$ to prove the result for any solution $\psi\in L^2(\R^N)$ (or even $\psi\in H^s(\R^N)$). (In fact the same issue arises in the local case $s=1$ discussed in the introduction.) We note that restricting to a finite number of $l$ values is reasonable when one has applications to atomic physics in mind.
    \e{enumerate}
\e{rmk}

\paragraph{Proof Strategy.}
The first two steps are analogues of the two steps in the proof of Theorem \ref{thm:Holder1}, but they are modified to exploit the spherical symmetry. We achieve this in part by invoking properties of hypergeometric functions, though more elementary approaches should also be possible (but will likely be more cumbersome).
\be{itemize}
    \item In step one, we suppose that $\psi$ is radial and run an analogous bootstrap argument as in the proof of Theorem \ref{thm:Holder1}. The localization-on-the-right Lemma \ref{lm:pwiserad} requires more work than its analogue, Lemma \ref{lm:pwise}.
		\item In step two, we still assume that $\psi$ is radial and conclude H\"older continuity from boundedness via the smoothing properties of the Riesz kernel. 
    \item To conclude, we use the Frank-Lenzmann observation mentioned above to remove the assumption that $\psi$ is radial (i.e.\ we ``trade off'' $l$ for dimension). 
\e{itemize}


\subsection{Further results}
We discuss two further results. First, Theorem \ref{thm:Holder1local} concerns the natural generalization of Theorem \ref{thm:Holder1} to finite-measure domains $\Om\subset\R^N$. The proof uses a non local version of the Leibniz rule and the main work is to control the ``localization error'' that it produces. 




\be{thm}[Finite-measure domains]
\label{thm:Holder1local}
Let $\Om\subset\set{R}^N$ be open and of finite measure. Let $s<1$ and $N,B,B',\al$ be as in Theorem \ref{thm:Holder1} with $B\cc \Om$. Let $f,V:\R^N\to \C$ with $f \in L^p(B)\cap L_{\mathrm{loc}}^1(\Om)$ and $V\in L^p(B)$ for some $1< p\leq \it$. Let $\psi$ be a distributional solution to
\beq
\label{eq:FSEOm}
		(-\Del)^s \psi + V\psi=f\quad \text{on } \Om,
\eeq
in the sense of Definition \ref{defn:distributional}. For $1/2\leq s<1$, assume additionally that $\psi\in W^{s',1+\de}(\Om)$ for some $s'>2s-1$ and $\de>0$. Then, the same conclusions as in Theorem \ref{thm:Holder1} hold.
\e{thm}

\textbf{Proof Strategy.}
\be{itemize}
    \item First, we extend $\psi$ to the function $\zeta \psi$ defined on $\set{R}^N$, where $\zeta\in C_0^\it(\R^N)$ is a cutoff function with
    \beq
    \label{eq:zetasupp}
        B\cc\supp\,\zeta\cc \Om
    \eeq
    and $\zeta\equiv 1$ on $B$.
The idea is to apply Theorem \ref{thm:Holder1} to $\zeta\psi$. This requires finding the equation satisfied by $\zeta\psi$ on $\set{R}^N$ (see \eqref{eq:FSEwhole}). (For the ordinary Laplacian, the Leibniz rule gives $-\De(\zeta\psi)= \zeta (-\De)\psi+\curly{E}$ with the ``localization error'' $\curly{E}=2\nabla\psi \nabla \zeta+\psi\De \zeta$.) We give a non local variant of the Leibniz rule, Proposition \ref{prop:nlleibniz}, which produces a different ``localization error''.

    \item Lemma \ref{lm:locerror} controls the localization error. The upshot is that $(-\De)^s(\zeta\psi)$ is locally bounded, and when $1/2\leq s<1$ its local integrability is as good as that of $(-\De)^s\psi$. This allows us to run a modified bootstrap argument. 
\e{itemize}

\be{rmk} 
\be{enumerate}[label=(\roman*)]
\item
For the definition of the fractional Sobolev space $W^{s',1+\de}(\Om)$, see e.g.\ \cite{Hitchhiker}. We need this additional assumption for $1/2\leq s<1$ to have an a priori estimate on the localization error (to understand this intuitively, note that the localization error $\curly{E}$ in the case $s=1$ from above depends on $\nabla\psi$). Also, we remark that $0\leq 2s-1<s$ for $s<1$, so the additional assumption is satisfied e.g.\ for any $\psi\in H^s(\R^N)$.
\item Note that we exclude the case $s=1$ in Theorem \ref{thm:Holder1local}. This case is technically easier, because the non local Leibniz rule can be replaced by the standard Leibniz rule, but the argument is different in that case (one cannot use \eqref{eq:nlleibniz2}, not even with a principal value). Since our main interest is in the non local situation, we do not consider $s=1$ here.
\e{enumerate}
\e{rmk}

The second additional result is that the H\"older exponents derived in Theorems \ref{thm:Holder1} and \ref{thm:Holder1local} are optimal.

\be{prop}[Optimality]
\label{thm:optimality}
The H\"older exponents $\al$, resp.\ $\al-1$ derived in Theorems \ref{thm:Holder1} and \ref{thm:Holder1local}
are \emph{optimal} for given values of the parameters $s,N,p$.

For example, let $\eps>0$ and let $N,s,p,B,B'$ be as in Theorem \ref{thm:Holder1} (I) with $\al=2s-N/p<1$. Then, there exist $f\in L^p(B)\cap \curly{L}_{N/2s}(\R^N)$ and a compactly supported distributional
solution $\psi$ to $(-\Del)^s \psi=f$ on $\set{R}^N$ such that $\psi\notin \tilde C^{0,\al+\eps}(B')$.
Analogous results hold in the other cases.
\e{prop}

\paragraph{Proof strategy.}
 We have $(-\Del)^s |x|^{\beta} = C(\beta) |x|^{\beta-2s}$ for any $0<\beta<2s$, see Proposition \ref{prop:scaling}. The argument is then roughly as follows: Let $\eps>0$. On the one hand, the function $|x|^{2s-N/p+\eps}$
 fails to be H\"older continuous of order $2s-N/p+2\eps$ at the origin. On the other hand, by the relation above, it corresponds to an $f$ of the form $|x|^{-N/p+\eps}$, which is in $L^p$ on any ball.

To make this heuristic rigorous, we multiply $|x|^\beta$ by a cutoff function $\zeta\in C_0^\it(\R^N)$ to get a compactly supported $\psi$. That $(-\De)^s(\zeta|x|^\beta)$ has the same local integrability as $|x|^{\beta-2s}$ is part of the proof of Theorem \ref{thm:Holder1local} as discussed above, see also Appendix \ref{app:optimality}.

\be{rmk}
Note that the functions $|x|^{2s-N/p+\eps/2}$ fail to have better H\"older regularity \emph{only at the origin}. This explains why Proposition \ref{thm:optimality} does not contradict Theorem \ref{thm:radial}.
\e{rmk}

\section{Proof of Theorem \ref{thm:Holder1}}

\subsection{Localizing on the right}
The proof is based on a bootstrap procedure which gives successively better $L^q$ regularity of $\psi$ on a sequence of balls sitting
 between $B'$ and $B$ until eventually $q=\it$, see Lemma \ref{lm:step1} below. To run the bootstrap argument, we need to localize $\psi$ appropriately at each step and for this we use the following lemma, which is based on an idea of Silvestre \cite{Silvestre07}.


\be{lm}
\label{lm:pwise}
Let $\psi$ be a distributional solution in the sense of Definition \ref{defn:distributional} to
 \beqs
		(-\Del)^s \psi + V\psi=f \quad \text{on } \set{R}^N.
\eeqs
Let $B_{\mathrm{small}}, B_{\mathrm{large}}$ be open balls such that $
 B_{\mathrm{small}} \cc B_{\mathrm{large}}$. Let $\eta\in C_0^\it(\R^N)$ be a cutoff function with $\eta\equiv 1$ on $B_{\mathrm{large}}$.
  Define $\psi_{\mathrm{loc}}$ by
  \beq
   \label{eq:psi1defn}
    \psi_{\mathrm{loc}}= k_{2s} *(g\eta),\qquad g=-V\psi+f.
  \eeq
Suppose that $g\in \curly{L}_{N/(2s)}(\R^N)$. Then, for any multi-index $\beta$, it holds that $D^\beta(\psi-\psi_{\mathrm{loc}}) \in L^\it(B_{\mathrm{small}})$.
\e{lm}

The upshot is that studying the local regularity of $\psi$ is the same as studying the local regularity of
$\psi_{\mathrm{loc}}$. The latter only depends on the values that $\psi$ takes on $\supp\,\eta$ and this is important for the bootstrap argument to work.


\be{proof}
 By Lemma \ref{lm:riesz}, $\psi_{\mathrm{loc}}$ satisfies
 \beqs
    (-\De)^s \psi_{\mathrm{loc}}=g\eta \quad \text{on } \set{R}^N.
 \eeqs
 in the distributional sense. By linearity, our assumption on $\psi$ gives
 \beqs
    (-\De)^s (\psi-\psi_{\mathrm{loc}}) = g(1-\eta)\quad \text{on } \set{R}^N.
 \eeqs
 in the distributional sense. By our assumption on $g$, Proposition \ref{prop:riesz} yields $w,w_1,\ldots,w_N\in\set{C}$ such that
 \beqs
    \psi(x)-\psi_{\mathrm{loc}}(x)= w+\sum_{j=1}^N w_j x_j +k_{2s} *  (g (1-\eta))(x).
 \eeqs
 The first two terms clearly have bounded derivatives of all orders everywhere, so it suffices to consider the third term. Let $x\in B_{\mathrm{small}}$. Note that $1-\eta\equiv 0$ on $B_{\mathrm{large}}$ and $0\leq \eta\leq 1$ everywhere. If either $N\neq 2s$ or $\beta\neq (0,\ldots,0)$, we have
\beq
\label{eq:above}
		|D^\beta k_{2s}*(g (1-\eta))(x)|
\leq C\int_{\set{R}^N\setminus B_{\mathrm{large}}} \frac{1}{|x-y|^{N-2s+|\beta|}} |g(y)|\dy.
\eeq
We make the general observation that $|x|\leq C$ and $|x-y|\geq C'$ implies that $|x-y|\geq C''\jap{y}$ with $C''$ depending only on $C'$ and $C''$. This yields
$$
\int_{\set{R}^N\setminus B_{\mathrm{large}}} \frac{1}{|x-y|^{N-2s+|\beta|}} |g(y)|\dy
\leq C\int_{\set{R}^N\setminus B_{\mathrm{large}}} \jap{y}^{2s-N-|\beta|} |g(y)|\dy.
$$
By H\"older's inequality, the bound is finite for any $g\in \curly{L}_{N/(2s)}(\set{R}^N)$ and this proves the claim in that case. 

It remains to consider the case $N=2s$ and $\beta=(0,\ldots,0)$, in which we estimate instead 
$$
C|\log|x-y||+C'\leq C''\max\{|x-y|^\de,|x-y|^{-\de}\}\leq C''' |y|^\de\leq C'''\jap{y}^\de,
$$
for any $\de>0$. The claim then follows from Definition \ref{defn:weightedL1} of $\curly{L}_1(\R^N)$.
\e{proof}

\subsection{Proof of boundedness by a bootstrap argument}

\be{lm}
\label{lm:step1}
Under the assumptions of Theorem \ref{thm:Holder1}, let $B_\IT$ be a ball such that $B'\cc B_\IT\cc B$. Then, $\psi\in L^\it(B_\IT)$.
\e{lm}

\be{proof}
 Let $B_1,B_{1/2}$ be balls such that $B_\IT\cc B_1 \cc B_{1/2}\cc B$. Let $\zeta_1\in C_0^\it(\R^N)$
 be a cutoff function supported in $B$ with $\zeta_1\equiv 1$ on $B_{1/2}$. Define
 \beq
 \label{eq:psi1}
  \psi_1: = k_{2s}* (g\zeta_1),\qquad \eps:=2s/N-p^{-1}.
 \eeq
 Note that $\eps>0$ in both cases (I),(II) and this is what drives the bootstrap procedure. By the definition of a distributional
 solution, $\psi\in L^{q_0}(B)$ with $q_0=1$. (We use $q_0$ instead of $1$ to elucidate the structure of the bootstrap argument.)
By H\"older's inequality and $g=-V\psi+f$ we have $g \zeta_1\in L^r(B)$ with $r^{-1}=p^{-1}+q_0^{-1}$.
Let $N>2s$. Applying Young's inequality to $\psi_1$ from \eqref{eq:psi1} yields $\psi_1\in L^{q_1}(B)\subset L^{q_1}(B_1)$ for all $q_1\geq 1$ satisfying
\beqs
		q_1^{-1} > -1 +\frac{N-2s}{N} + r^{-1} = q_0^{-1} - \eps.
\eeqs
and we pick $q_1^{-1}=q_0^{-1}-\eps/2$ for definiteness. 

By assumption, we have $g\in \curly{L}_{N/(2s)}(\set{R}^N)$. Hence, we can apply Lemma \ref{lm:pwise} with the choices $B_{\mathrm{small}}=B_1, B_{\mathrm{large}}=B_{1/2}$ and $\eta=\zeta_1$,
to get that $\psi-\psi_1$ is uniformly bounded on $B_1$. Thus, $\psi_1\in L^{q_1}(B_1)$ implies $\psi\in L^{q_1}(B_1)$ as well.

 If $q_1^{-1}\geq \eps$, we repeat the procedure as follows: We take open balls $B_2,B_{3/2}$ such that $B_\IT\cc B_2\cc B_{3/2}\cc B_1$ and a cutoff function $\zeta_2\in C_0^\it(\R^N)$ supported in $B_1$ with $\zeta_2\equiv 1$ on $B_{3/2}$. We define
 \beqs
  \psi_2: = k_{2s}* (g\zeta_2).
 \eeqs
By combining H\"older's inequality and Young's inequality, we get $\psi_2\in L^{q_2}(B_2)$ for all $q_2\geq 1$ satisfying $q_2^{-1}>q_1^{-1}-\eps$ and we pick $q_2^{-1}= q_1^{-1}-\eps/2$ for definiteness. As before, Lemma \ref{lm:pwise} implies that $\psi-\psi_2$ is uniformly bounded on $B_2$ and consequently $\psi\in L^{q_2}(B_2)$. After finitely many (say $k$) iterations, we obtain $\psi\in L^{q_{k}}(B_{k})$ with $q^{-1}_{k}<\eps$ with $B_\it\cc B_k\cc B$. 

In the case $N=2s$, the argument can be modified as follows. For $x,y$ in some fixed compact sets, we have
$$
C|\log|x-y||+C'\leq C''|x-y|^{-\de}+C'''
$$
for any $\de>0$. If we choose $\de$ sufficiently small compared to $\eps=2s/N-p^{-1}>0$, Young's inequality gives a finite improvement of integrability on the $L^q$ scale as above. Consequently, after some number of iterations (called $k$), we have the same conclusion as in the case $N>2s$.\\

 To conclude that $\psi\in L^\it(B_\IT)$, we take a ball $B_\IT'$ such that $B_\IT\cc B_\IT' \cc B_k$ and a
cutoff function $\zeta_\it\in C_0^\it(\R^N)$ supported in $B_k$ with $\zeta_\it\equiv 1$ on $B_\IT'$. We define
\beqs
\psi_\it := k_{2s}*(g\zeta_\it).
\eeqs
By H\"older's inequality and the fact
that $\zeta_\it$ is supported in $B_k$, we have
\begin{align*}
        \sup_{x\in B_\IT'} |\psi_\it(x)|&
        \leq C \sup_{x\in B_\IT'} \int |k_{2s}(x-y)| |g(y) \zeta_\it(y)| \dy\\
          &\leq \|g\zeta_\it\|_{L^r(B_{k})} \|k_{2s}\|_{L^{r'}(2B_{k})}.
\end{align*}
The first term is finite (by H\"older's inequality) if we choose $r^{-1}=p^{-1}+q_{k-1}^{-1}$. Since $q^{-1}_{k}<\eps$, we have $r^{-1}<2s/N$, which ensures that $r'(N-2s)<N$ and hence that $\|k_{2s}\|_{L^{r'}(2B_{k-1})}$ is also finite. This includes the case $N=2s$, in which we use that $p^{-1}=1-\eps$ gives $r>1$ and therefore $r'<\it$ (and any finite power of the logarithm is locally integrable). This proves $\psi_\it\in L^\it(B_\IT')\subset L^\it(B_\IT)$.  By another application of Lemma \ref{lm:pwise}, we get $\psi\in L^\it(B_\IT)$ and Lemma \ref{lm:step1} follows.
 \e{proof}

\subsection{Concluding the continuity from boundedness}
\label{ssect:conclusion}
We apply Lemma \ref{lm:step1} to obtain that $\psi\in L^\it(B_\IT)$ for some ball $B_\IT$ with $B'\cc B_\IT\cc B$. We will now use the smoothing properties of the Riesz kernel to conclude the claimed H\"older/almost Lipschitz continuity of $\psi$ on $B'$. 


\be{proof}

Let $p>N/(2s)$. We localize on the right again. That is, we take a ball $B_H$
with $B'\cc B_{H} \cc B_\it$ and a cutoff function $\zeta_{H}\in C_0^\it(\R^N)$ supported on $B_\IT$ with $\zeta_H\equiv 1$ on $B_{H}$. We define
 \beq
 \label{eq:proof12}
 		\psi_{H} := k_{2s}*(g \zeta_{H})
 \eeq
 where $g=-V\psi+f$. We will prove that $\psi_{H}$ has the claimed H\"older/almost-Lipschitz continuity on $B'$.
  By Lemma \ref{lm:pwise}, $\nabla(\psi-\psi_H)$ is uniformly bounded on $B'$ and so $\psi-\psi_H$ is uniformly Lipschitz continuous on $B'$. Hence, $\psi$ will have the same H\"older/almost-Lipschitz continuity as $\psi_H$ on $B'$.

\dashuline{Proof of statements (I),(II.A).} Let $x_1,x_2\in B'$ and recall that $\supp\,\zeta_H\cc B_\IT$.
By applying H\"older's inequality twice, we have
 \beq
 \begin{aligned}
 \label{eq:proof14}
       &|\psi_{H}(x_1)-\psi_{H}(x_2)|\\
       \leq&  C \|g\zeta_{H}\|_{L^p(B_\IT)} \l(\int_{B_\IT} \l| k_{2s}(x_1-y)-k_{2s}(x_2-y) \r|^{p'}\dy\r)^{1/p'}\\
       \leq& C(\|V\|_{L^p(B)}\|\psi\|_{L^\it(B_\it)}+\|f\|_{L^p(B)}) \l(\int_{B_\IT} \l| k_{2s}(x_1-y)-k_{2s}(x_2-y) \r|^{p'}\dy\r)^{1/p'}.
  \end{aligned}
  \eeq
Note that $p>1$ and so $p'<\it$. We introduce the unit vector $e=(x_1-x_2)/|x_1-x_2|$ and let $R>2$ be such that $B_H+B_\IT\cc R \curly{B}$, where $\curly{B}$ denotes the unit ball centered at the origin. By scaling and translation, we have
 \begin{align}
 \nonumber
      &\l(\int_{B_\it} \l| k_{2s}(x_1-y)-k_{2s}(x_2-y) \r|^{p'}\dy\r)^{1/p'}\\
  \label{eq:firstfactor}
      &\leq  |x_1-x_2|^{N/p'-(N-2s)} \l(\int_{|x_1-x_2|^{-1}R \curly{B}} \l| k_{2s}(e-y)-k_{2s}(y) \r|^{p'}\dy\r)^{1/p'}
 \end{align}
 We emphasize that \eqref{eq:firstfactor} also holds in the case $N=2s$ (using $\log a-\log b=\log(a/b)$). We observe that $N/p'-(N-2s)=2s-N/p=\al$ is already the claimed H\"older exponent. To prove statement (I), it remains to control the second
 factor in \eqref{eq:firstfactor}, either by a constant, when $\al<1$, or by a constant times $\log(|x_1-x_2|^{-1})$, when $\al=1$.\\

 Since $\psi_H$ is bounded on $B'$, we may assume without loss of generality that $|x_1-x_2|\leq 1$. The crucial idea is to cut the integral as follows
 \beq
 \label{eq:crucialcut}
 \begin{aligned}
    &\int_{|x_1-x_2|^{-1}R \curly{B}} \l| k_{2s}(e-y)-k_{2s}(y) \r|^{p'}\dy\\
   =&\int_{|x_1-x_2|^{-1}R \curly{B}\setminus 2 \curly{B}} \l| k_{2s}(e-y)-k_{2s}(y)  \r|^{p'}\dy+ \int_{2\curly{B}} \l| k_{2s}(e-y)-k_{2s}(y)  \r|^{p'}\dy
 \end{aligned}
 \eeq
We observe that the second integral is finite. Indeed, when $N=2s$, this holds because $p'<\it$ and any power of a logarithm is locally integrable.
When $2s<N$, finiteness follows from the elementary inequality $(a+b)^{p'}\leq 2^{p'-1}(a^{p'}+b^{p'})$ for any $a,b>0$ and the fact that
  $(N-2s)p'<N$, which is equivalent to our assumption $p>N/(2s)$.
	
	Consider the first integral in \eqref{eq:crucialcut}. Since the Riesz kernel is differentiable away from its singularity, we can apply the mean value theorem to get
\beq
\label{eq:rho}
\begin{aligned}
&\int_{|x_1-x_2|^{-1}R \curly{B}\setminus 2 \curly{B}}
\l| k_{2s}(e-y)-k_{2s}(y) \r|^{p'}\dy\\
\leq& C\int_{|x_1-x_2|^{-1}R \curly{B}\setminus 2 \curly{B}}
\l| \frac{1}{(|y|-1)^{(N-2s+1)p'}} \r|^{p'}\dy
\leq    C \int_{|x_1-x_2|^{-1}R \curly{B}\setminus 2 \curly{B}} \frac{1}{|y|^{(N-2s+1)p'}}\dy\\
=& C \int_2^{R|x_1-x_2|^{-1}} \rho^{N-1-(N-2s+1)p'}  \d\rho\\
=& C+
    \be{cases}
        C' |x_1-x_2|^{(N-2s+1)p'-N},\quad \text{if } 2s-N/p<1\\
        C' \log(|x_1-x_2|^{-1}),\quad \text{if } 2s-N/p=1.
    \e{cases}
\end{aligned}
\eeq
Observe that $2s-N/p<1$ gives $(N-2s+1)p'-N>0$. This proves statements (I) and (II.A).\\

\dashuline{ Proof of statement (II.B).}
Let $1/2<s\leq 1$ and $p>N/(2s-1)$. Let $B_H,\zeta_H,\psi_H$ be as before. We first show that $\psi_H\in C^1(B_\IT)$ and identify its derivative, see \eqref{eq:proof13}. For any unit vector $e$ and any $x\neq0$, we have
\beq
\label{eq:differencequotientlimit}
    \lim_{h\rightarrow0} \frac{k_{2s}(x+he)-k_{2s}(x)}{h} = C_0\frac{\scp{e}{x}}{|x|^{N-2s+2}}
\eeq
where the value of $C_0$ follows from \eqref{eq:rieszdefn}. Hence, by the mean value theorem,
\beq
\label{eq:differencequotientcontrol}
\begin{aligned}
\l| \frac{k_{2s}(x+he)-k_{2s}(x)}{h}\r| \leq  C \frac{1}{||x|-1|^{N-2s+1}}
\end{aligned}
\eeq
for all $|h|\leq 1$. Recall that $\psi_H(x) =k_{2s}*(g\zeta_H)(x)$. Let $x\in B_\IT$ and consider
\begin{align*}
    \frac{k_{2s}*(g\zeta_H)(x+he)-k_{2s}*(g\zeta_H)(x)}{h}
    = \int  \frac{k_{2s}(x-y+he)-k_{2s}(x-y)}{h} g(y)\zeta_1(y)\d y.
\end{align*}
We have $g\zeta_H\in L^p(B_\IT)$, by H\"older's inequality and $\psi\in L^\it(B_\IT)$. Using \eqref{eq:differencequotientcontrol}, $\supp\,\zeta_{H}\cc B_\IT$
and the fact that $p>N/(2s-1)$, we see that the integrand is bounded by an integrable function uniformly in $h$ for all $0<|h|\leq 1$. Thus,
dominated convergence and \eqref{eq:differencequotientlimit} imply that $\psi_{H}\in C^1(B_\IT)$ with
\beq
\label{eq:proof13}
    \nabla \psi_{H}(x) = C_0\int \frac{x-y}{|x-y|^{N-2s+2}} g(y)\zeta_{H}(y)\dy.
\eeq
One then follows the same line of argumentation as for statements (I) and (II.A). The fact that \eqref{eq:proof13} is vector valued only requires modifying the proof in the first step of \eqref{eq:rho}. One uses the following estimate which holds for all $y\in\R^N$ with $|y|>2$ and all unit vectors $e$ (concretely $e=\frac{x_1-x_2}{|x_1-x_2|}$):
$$
\begin{aligned}
\l|\frac{y-e}{|y-e|^{N-2s+2}}-\frac{y}{|y|^{N-2s+2}}\r|
\leq &\frac{1}{(|y|-1)^{N-2s+2}}+|y|\l|\frac{1}{|y-e|^{N-2s+2}}-\frac{1}{|y|^{N-2s+2}}\r|\\
\leq &C\frac{1}{|y|^{N-2s+2}}.
\end{aligned}
$$
In the second step, we used the mean value theorem and $|y|-1>|y|/2$. The remaining details in the proof of statement (II.B) are left to the reader.
\e{proof}



\section{Proof of Theorem \ref{thm:radial}}
Let $N\geq 2$. Since the proof is simpler in the case $s=1$ and was sketched in the introduction, we assume $2s<N$ henceforth.

Consider $k_{2s}*g$ with $g$ radial. Independently of $g$, the angular integral in this expression reads
\beq
\label{eq:Phidefn}
    \Phi(t):=|\S^{N-2}| \int_{-1}^1 \frac{(1-u^2)^{\frac{N-3}{2}  }}{(1+t^2-2 t u)^{N/2-s}}  \d u,
\eeq
with the convention that $|\S^0|=2$. Here $t$ is a multiple of the radial integration variable, see \eqref{eq:precise} for the precise expression.

 Now observe that $\Phi$ is only singular when the denominator
\beqs
1+t^2-2 t u=(1-tu)^2+t^2(1-u^2)
\eeqs
vanishes, i.e.\ when $t=u=1$. While at first sight the singularity may appear to be of order $N-2s$ (we are considering the Riesz kernel $k_{2s}$ after all), the point $u=1$ is ``dampened'' by the other factor $(1-u^2)^{\frac{N-3}{2}}$ in \eqref{eq:Phidefn}. It turns out that this reduces the order of the singularity to $1-2s$, as can be seen in \eqref{eq:HG}, which shows that the \emph{effective dimension is indeed one}.


\subsection{Properties of $\Phi$}

For odd values of $N$, the above argument can be made rigorous by repeated integration by parts in $u$ in \eqref{eq:Phidefn}. To include even values of $N$, one can either use elementary estimates or, as we do here, invoke known properties of hypergeometric functions. This is inspired by the proof of Lemma 3.1 in \cite{FrankSeiringer08gs}.

\be{prop}
\label{prop:Phi}
\be{enumerate}[label=(\roman*)]
\item $\Phi:(-1,1)\to \R_+$ is analytic. 
\item For all $t\in (0,1)$,
\beq
\begin{aligned}
\label{eq:HG}
\Phi(t)=\phi_1(t)+K_{2s}(t)\phi_2(t),\qquad \text{with }
K_{2s}(t):=
\be{cases}
|t|^{2s-1},\quad &\text{if } s\neq 1/2\\
\log|t|,\quad &\text{if } s=1/2,
\e{cases}
\end{aligned}
\eeq
where $\phi_1,\phi_2:(0,\sqrt{2})\to \R$ are smooth functions.
\item For all $t>0$, we have the functional equation $\Phi(t^{-1})=t^{N-2s}\Phi(t)$.
\e{enumerate}
\e{prop} 

\be{rmk}
\be{enumerate}[label=(\roman*)]
\item
$K_{2s}$ is (up to constants) the Riesz kernel $k_{2s}$ with $N=1$.
\item
We emphasize that $\phi_1,\phi_2$ are smooth on the larger interval $(0,\sqrt{2})$, in other words \eqref{eq:HG} isolates the potentially singular part $K_{2s}(t)$ of $\Phi(t)$. 
\item We will generally use the analyticity/smoothness of the functions appearing in (i) and (ii) only to conclude that they and their derivatives are \emph{bounded} on any compact subset of their domain of definition.
\e{enumerate}
\e{rmk}

\be{proof}
We derive the properties of $\Phi$ through its connection to $F(a,b;c;z)$, the ordinary hypergeometric function (also denoted ${}_2 F_1$). By formula (3.665) in \cite{GradshteynRyzhik} we have
\beq
	\label{eq:hyper1}
    \Phi(t) = \l|\S^{N-2}\r| B\l(\frac{N-1}{2},\frac{1}{2}\r)F\l(\frac{N}{2}-s,1-s;\frac{N}{2};t^2\r),\qquad t\in(-1,1),
    \eeq
where $B(a,b)$ denotes the Beta function. The hypergeometric function $F(a,b;c;z)$ has a series representation around $z=0$ with a radius of convergence equal to one, see \cite{AbramowitzStegun} (15.1.1). The map $(-1,1)\to (-1,1)$, $t\mapsto t^2$ is analytic, so that $\Phi:(-1,1)\to \R_+$ is a composition of analytic functions and hence itself analytic.

We come to statement (ii). Let $s\neq 1/2$ first. By formula (15.3.6) in \cite{AbramowitzStegun} and the fact that $t\mapsto F(a,b;c;1-t^2)$ is analytic on $(0,\sqrt{2})$, there exist (explicit) analytic functions $\phi_1,\tilde\phi_2:(0,\sqrt{2})\to \R$ such that
$$
\Phi(t)=\phi_1(t)+(1-t^2)^{2s-1}\tilde\phi_2(t),\qquad t\in(-1,1).
$$
Note that $(1-t^2)^{2s-1}=(1-t)^{2s-1}(1+t)^{2s-1}$. Since the map $t\mapsto (1+t)^{2s-1}$ is smooth on $(0,\sqrt{2})$, we can set $\phi_2(t)=(1+t)^{2s-1}\tilde\phi_2(t)$ to get \eqref{eq:HG} with smooth $\phi_1,\phi_2$. 

For $s=1/2$, one is in a degenerate case (the hypergeometric function in \eqref{eq:hyper1} is of the form $F(a,b;a+b;z)$).  By formula (15.3.10) in \cite{AbramowitzStegun} we have in this case
$$
\Phi(t)=\sum_{n=0}^\it (A_n+B_n\log(1-t^2))(1-t^2)^n=:\tilde \phi_1(t) +\log(1-t^2)\phi_2(t),\qquad t\in(-1,1),
$$
for some explicit coefficients $A_n,B_n$. They are such that the two series above converges absolutely for all $t\in (0,\sqrt{2})$, which in turn defines the analytic functions $\tilde\phi_1,\phi_2:(0,\sqrt{2})\to \R$ in the second equality. Note that $\log(1-t^2)=\log(1-t)+\log(1+t)$. Since the map $t\mapsto \log(1+t)$ is smooth on $(0,\sqrt{2})$, we can set $\phi_1(t)=\log(1+t)+\tilde\phi_1(t)$ to get \eqref{eq:HG} when $s=1/2$. This proves statement (ii).

Finally, statement (iii) follows directly from the definition \eqref{eq:Phidefn} of $\Phi$.
\e{proof}

\subsection{Localizing on the right}
The following is a radial analogue of Lemma \ref{lm:pwise}. Recall that for a radial function $F$ we regularly abuse notation by identifying $F(x)\equiv F(|x|)$.

\be{lm}
\label{lm:pwiserad}
Let $V,f:\R_+\to \C$ be radial. Let $\psi$ be a distributional solution to
\beqs
	(-\De)^s\psi + V\psi =f\quad \text{on }\set{R}^N
\eeqs
and suppose that $\psi$ is radial. Consider two nested intervals $[a_{\mathrm{small}},b_{\mathrm{small}}]\cc [a_{\mathrm{large}},b_{\mathrm{large}}]$ with $a_{\mathrm{large}}>0$.
Let $\eta\in C_0^\it(\R_+)$ be a cutoff function with $\eta\equiv 1$ on
$[a_{\mathrm{large}},b_{\mathrm{large}}]$. Define
\beqs
	\psi_{\mathrm{loc}} =k_{2s}*(g \eta)\quad \text{on }\set{R}^N
\eeqs
where $g=-V\psi+f$. Suppose that $g$ satisfies \eqref{eq:reg} with $l=0$, i.e.\
 \beq
 \label{eq:regg}
 g\in 
 \be{cases}
 \curly{L}_{1/(2s)}(\R_+; r^{N-1}\d r),\quad &\text{if } 0<s\leq 1/2\\
 \jap{\cdot}^{N-2s}L^1(\R_+; r^{N-1}\d r),\quad &\text{if } 1/2<s\leq 1.
 \e{cases}
 \eeq
Then $\psi-\psi_{\mathrm{loc}}$, $(\psi-\psi_{\mathrm{loc}})'$ and $(\psi-\psi_{\mathrm{loc}})''$ are bounded on $[a_{\mathrm{small}},b_{\mathrm{small}}]$.
\e{lm}

 \be{proof}
 Let $\rho\in [a_{\mathrm{small}},b_{\mathrm{small}}]$. By changing to spherical coordinates, we have
 \beq
 \label{eq:precise}
\begin{aligned}
 \psi(\rho)-\psi_{\mathrm{loc}}(\rho)=k_{2s}*\l(g(1-\eta) \r)(\rho)= \rho^{2s-N} \int_0^\it  \Phi\l(\frac{r}{\rho}\r) g(r)(1-\eta(r)) r^{N-1}\d r,
\end{aligned}
 \eeq
with $\Phi$ as in \eqref{eq:Phidefn}. Using that $1-\eta\equiv 0$ on $[a_{\mathrm{large}},b_{\mathrm{large}}]$ and $0\leq \eta\leq 1$ everywhere, we get
\beq
\begin{aligned}
\label{eq:locrad1}
&|\psi(\rho)-\psi_{\mathrm{loc}}(\rho)|\\
&\leq \rho^{2s-N}\int_0^{a_{\mathrm{large}}}  \Phi\l(\frac{r}{\rho}\r)|g(r)|r^{N-1}\d r + \rho^{2s-N} \int_{b_{\mathrm{large}}}^{\it}  \Phi\l(\frac{r}{\rho}\r) |g(r)|r^{N-1}\d r.
\end{aligned}
\eeq
Consider the first integral. Observe that the argument of $\Phi$ satisfies 
$$
\frac{r}{\rho}\in \l[0,\frac{a_{\mathrm{large}}}{\rho}\r].
$$
Recall $\rho\in [a_{\mathrm{small}},b_{\mathrm{small}}]$. Our assumption $[a_{\mathrm{small}},b_{\mathrm{small}}]\cc [a_{\mathrm{large}},b_{\mathrm{large}}]$, implies that $a_{\mathrm{small}}\geq a_{\mathrm{large}}>0$. Hence $a_{\mathrm{large}}/\rho\in (0,1)$ and Proposition \ref{prop:Phi} (i) implies that $\Phi$ is bounded on $[0,a_{\mathrm{large}}/\rho]$. Moreover, $\rho^{2s-N}\leq C$ since $\rho> a_{\mathrm{large}}>0$. Therefore
\beq
\label{eq:firstint}
\rho^{2s-N}\int_0^{a_{\mathrm{large}}} \Phi\l(\frac{r}{\rho}\r)|g(r)|r^{N-1}\d r\leq C \int_0^{a_{\mathrm{large}}}|g(r)| r^{N-1}\d r.
\eeq
Consider the second integral in \eqref{eq:locrad1}. We will relate it back to the previous case by using the functional equation from Proposition \ref{prop:Phi} (iii), i.e.
\beq
\label{eq:Phiscaling}
	\Phi\l(\frac{r}{\rho}\r) = \l(\frac{\rho}{r}\r)^{N-2s} \Phi\l(\frac{\rho}{r}\r).
\eeq
It allows us to repeat the arguments that gave \eqref{eq:firstint} for the second integral in \eqref{eq:locrad1}. The result is
\beqs
    |\psi(\rho)-\psi_{\mathrm{loc}}(\rho)|
    \leq C \int_0^{a_{\mathrm{large}}}|g(r)| r^{N-1}\d r + C' \int_{b_{\mathrm{large}}}^{\it} \jap{r}^{2s-N} |g(r)| r^{N-1}\d r.
\eeqs
Both integrals are finite by our assumption \eqref{eq:regg} on $g$. Since the constants on the right hand side do not depend on the choice of $\rho\in [a_{\mathrm{small}},b_{\mathrm{small}}]$, we get that $\psi-\psi_{\mathrm{loc}}$ is bounded on $[a_{\mathrm{small}},b_{\mathrm{small}}]$.\\

We come to $(\psi-\psi_{\mathrm{loc}})'$. Taking derivatives in \eqref{eq:precise}, we obtain
\beq
\label{eq:derivative}
(\psi-\psi_{\mathrm{loc}})'(\rho)=(2s-N)\rho^{-1}(\psi(\rho)-\psi_{\mathrm{loc}}(\rho))+\rho^{2s-N}\int_0^\it  \frac{\d}{\d \rho}\Phi\l(\frac{r}{\rho}\r) g(r)(1-\eta(r)) r^{N-1}\d r.
\eeq
We just proved that $\psi-\psi_{\mathrm{loc}}$ is bounded on $[a_{\mathrm{small}},b_{\mathrm{small}}]$ and $\rho^{-1}\leq C$ on the same interval. Hence, the first term in \eqref{eq:derivative} is bounded. Following the development in \eqref{eq:locrad1}, the second term in \eqref{eq:derivative} can be split into two parts. The first integral from $0$ to $a_{\mathrm{large}}$ can be bounded by the same arguments as before (instead of $\Phi\in L^\it([0,a_{\mathrm{large}}/\rho])$ one uses $\Phi'\in L^\it([0,a_{\mathrm{large}}/\rho])$ by Proposition \ref{prop:Phi} (i)). The more interesting term is the other one:
\beq
\label{eq:interesting}
\int_{b_{\mathrm{large}}}^\it \frac{\d}{\d \rho}\Phi\l(\frac{r}{\rho}\r)g(r) r^{N-1}\d r,
\eeq
where we have already dropped some powers of $\rho$ since they are uniformly bounded on $[a_{\mathrm{small}},b_{\mathrm{small}}]$. By the chain rule and Proposition \ref{prop:Phi} (iii)
we compute 
\beq
\label{eq:Phiderivative}
\frac{\d}{\d \rho}\Phi\l(\frac{r}{\rho}\r)=-\rho^{-2}r\Phi'\l(\frac{r}{\rho}\r)
 =-\rho^{N-2s-1} (2s-N)r^{2s-N}\Phi\l(\frac{\rho}{r}\r)+\rho^{N-2s}r^{2s-N-1}\Phi'\l(\frac{\rho}{r}\r).
\eeq
Observe that any ration $\rho/r$ appearing in \eqref{eq:derivative} satisfies $\rho/r\leq b_{\mathrm{small}}/b_{\mathrm{large}}<1$. Using again that $\Phi,\Phi'$ are bounded on $[0,\beta]$ for any $0<\beta<1$, we get
$$
\int_{b_{\mathrm{large}}}^\it \frac{\d}{\d \rho}\Phi\l(\frac{r}{\rho}\r)g(r) r^{N-1}\d r
\leq C \int_{b_{\mathrm{large}}}^\it \jap{r}^{2s-N}g(r)r^{N-1}\d r<\it
$$
This yields the claimed boundedness of $(\psi-\psi_{\mathrm{loc}})'$ on $[a_{\mathrm{small}},b_{\mathrm{small}}]$. For $(\psi-\psi_{\mathrm{loc}})''$, the argument is similar and left to the reader (note that taking another $\frac{\d}{\d \rho}$ derivative in \eqref{eq:Phiderivative} can only produce more inverse powers of $r$ and these ameliorate the decay). 
\e{proof}


%

\subsection{Proof of boundedness}

The following lemma is a radial analogue of Lemma \ref{lm:step1}. Recall that by our convention an ``annulus'' is always centered at the origin. 

\be{lm}
\label{lm:step1rad}
Under the assumptions of Theorem \ref{thm:radial}, let $A_\it$ be an open annulus such that $A'\cc A_\it \cc A$. Assume further that $\psi$ is radial. Then, with the usual abuse of notation for radial functions $\psi\in L^\it(A_\it)$.
\e{lm}

\be{proof}
Let $A_1,A_{1/2}$ be open annuli such that $A_\it\cc A_1 \cc A_{1/2}\cc A$. Let $\zeta_1\in C_0^\it(\R_+)$ be a radial cutoff function supported in $A$ with $\zeta_1\equiv 1$ on $A_{1/2}$. We define
\beqs
    \psi_1 := k_{2s}*(g\zeta_1).
\eeqs
We denote the minimal/maximal radii of the annuli $A_1$ and $A$ as follows
\beq
\begin{aligned}
\label{eq:achoices}
a_{1}&=\min\setof{|y|}{y\in A_1},\quad b_{1}=\max\setof{|y|}{y\in A_1}\\
a&=\min\setof{|y|}{y\in A},\quad b=\max\setof{|y|}{y\in A}.
\end{aligned}
\eeq
Let $\rho\in [a_{1},b_{1}]$. In spherical coordinates, we have, with $\Phi$ as in \eqref{eq:Phidefn},
\beq
\label{eq:rad0}
  \psi_{1}(\rho)  =\rho^{2s-N} \int_{a}^{b}  \Phi\l(\frac{r}{\rho}\r)  g(r)\zeta_1(r)  r^{N-1}\d r.
\eeq
Note that we can restrict to $(\rho,r)$ values in the set \beq
\label{eq:Cdefn}
\curly{C}:=[a_{1},b_{1}]\times [a,b].
\eeq
The assumption $\min A'>0$ in Theorem \ref{thm:radial} implies $a,a_1>0$. Hence, all maps of the form $(\rho,r)\mapsto \rho^{\beta}r^{\gam}$ with $\beta,\gam\in \R$ are smooth on $\curly{C}$ and need not be considered further.
By Proposition \ref{prop:Phi}, the function $\Phi\l(r/\rho\r)$ is smooth on $\curly{C}$ away from the point $r=\rho$ and it behaves like the one dimensional Riesz kernel $K_{2s}$ near its singularity. Indeed, combining Proposition \ref{prop:Phi} (ii) and (iii), (and letting $s\neq 1/2$ for definitess), we have
\beq
\label{eq:Phicases}
\begin{aligned}
	 \Phi\l(\frac{r}{\rho}\r)=
\be{cases}
    \phi_1\l(\frac{r}{\rho}\r)+\l|\rho-r\r|^{2s-1}\rho^{1-2s}\phi_2\l(\frac{r}{\rho}\r), \quad &\text{if } 0<r<\rho\\
    \l(\frac{r}{\rho}\r)^{2s-N}\phi_1\l(\frac{\rho}{r}\r)+\l|\rho-r\r|^{2s-1}\rho^{1-2s}\l(\frac{r}{\rho}\r)^{1-N}\phi_2\l(\frac{\rho}{r}\r), \quad &\text{if } r> \rho.
\e{cases}
\end{aligned}
\eeq
Proposition \ref{prop:Phi} also yields a similar formula for $s=1/2$, but with $\l|\rho-r\r|^{2s-1}$ replaced by $\log|\rho-r|$. Recall that $\phi_1,\phi_2:(0,\sqrt{2})\to \R$ are smooth functions on an open neighborhood of $\curly{C}$. Consequently, they are \emph{uniformly bounded} in \eqref{eq:Phicases} for all pairs $(\rho,r)\in\curly{C}$.

Recall that all powers of $\rho$ and $r$ are smooth on $\curly{C}$. Hence, \eqref{eq:Phicases} (and its $s=1/2$ analogue) can be simplified as follows. There exist smooth functions $\vp_i:\curly{C}\to\R$, $i=1,2,3,4$ such that
\beq
\label{eq:PhiK}
\Phi\l(\frac{r}{\rho}\r)=
\be{cases}\vp_1(\rho,r)+K_{2s}(\rho-r)\vp_2(\rho,r), \quad &\text{if } r<\rho\\
\vp_3(\rho,r)+K_{2s}(\rho-r)\vp_4(\rho,r), \quad &\text{if } r> \rho
\e{cases}
\eeq
holds for all $0<s<1$. The one dimensional Riesz kernel $K_{2s}$ was defined in \eqref{eq:HG}.

\dashuline{ Case $1/2<s<1$.} Boundedness is immediate. Indeed, for $s>1/2$, $$K_{2s}(\rho-r)=|\rho-r|^{2s-1}\leq \l(|\rho|+|r|\r)^{2s-1}\leq C.$$ 
Since the $\vp_i$ in \eqref{eq:PhiK} are uniformly bounded, we find $|\Phi(r/\rho)|\leq C$ almost everywhere in $\curly{C}$ (namely except at $r=\rho$). Since all powers of $r$ and $\rho$ are bounded on $\curly{C}$, \eqref{eq:rad0} now gives
$$
|\psi_1(\rho)|\leq  C \int_{a}^{b} g(r)r^{N-1}\d r=C'.
$$
Here we used that $g\in L^1_{\mathrm{loc}}(\R_+;r^{N-1}\d r)$ by assumption. This proves $\psi_1\in L^\it(A_1)$. We then apply Lemma \ref{lm:pwiserad} with $\eta\equiv \zeta_1$ and intervals $[a_{\mathrm{small}},b_{\mathrm{small}}]=[a_1,b_1]$ and $[a_{\mathrm{large}},b_{\mathrm{large}}]=[a,b]$. We conclude that $\psi\in L^\it(A_1)$.

\dashuline{ Case $0<s\leq 1/2$.} We split the integral in \eqref{eq:rad0} into the regions $\{r<\rho\}$ and $\{r>\rho\}$ (and we ignore the null set $\{r=\rho\}$ henceforth). On each region, we can write $\Phi$ via \eqref{eq:PhiK} in terms of the one dimensional Riesz kernel $K_{2s}$ (up to additional smooth functions). In other words, we are essentially in the $N=1$ case of Lemma \ref{lm:step1}, (up to additional smooth functions). The proof then follows exactly the line of argument in the proof of Lemma \ref{lm:step1}, where every application of Lemma \ref{lm:pwise} is replaced by one of Lemma \ref{lm:pwiserad}. We leave the details to the reader.

\e{proof}
%
%

\subsection{Conclusion for radial solutions}

\be{proof}[Proof of Theorem \ref{thm:radial} for radial $\psi$]
Let $\psi$ be radial, or equivalently, let $\psi\in\curly{H}_0$. Take an open annulus $A_\it$ such that $A'\cc A_\it \cc A$. By Lemma \ref{lm:step1rad}, $\psi\in L^\it(A_\it)$.\\

Let $A_H$ be an open annulus with $A'\cc A_H\cc A_\it$. We define
$$
\begin{aligned}
a_H=&\min \setof{|y|}{y\in A_H},\quad b_H=\max \setof{|y|}{y\in A_H},\\
a_\it=&\min \setof{|y|}{y\in A_\it},\quad b_\it=\max \setof{|y|}{y\in A_\it},
\end{aligned}
$$
so that $[a_H,b_H]\cc[a_\it,b_\it]$. Note that by our assumption in Theorem \ref{thm:radial}, $a_H,a_\it>0$. Let $\zeta_H\in C_0^\it(\R_+)$ be a cutoff function supported in $[a_\it,b_\it]$ with $\zeta_H\equiv 1 $ on $[a_H,b_H]$. Define
\beqs
    \psi_{H} := k_{2s}*(g\zeta_{H}),\quad g=-V\psi+f
\eeqs
By the assumption \eqref{eq:reg} (recall that $l=0$), we can apply Lemma \ref{lm:pwiserad} to get that $(\psi-\psi_H)'$ and $(\psi-\psi_H)''$ are bounded on $[a',b']$. In particular, $\psi-\psi_H$ and $(\psi-\psi_H)'$ are uniformly Lipschitz continuous on $[a',b']$. Hence, it suffices to establish the claims (I) and (II) for $\psi_H$ instead of $\psi$.

Let $\rho\in [a_H,b_H]$. In spherical coordinates, we have (compare \eqref{eq:rad0})
\beq
\label{eq:psiH}
  \psi_{H}(\rho)  =\int_{a_\it}^{b_\it}  \rho^{2s-N} \Phi\l(\frac{r}{\rho}\r)  g(r)\zeta_H(r)  r^{N-1}\d r.
\eeq

\dashuline{ Proof of statements (I),(II.A).}
Let $\rho_1,\rho_2\in [a_H,b_H]$. Without loss of generality, we may assume $\rho_1\leq \rho_2$. We want to bound $\psi_{H}(\rho_1)-\psi_{H}(\rho_2)$. Using \eqref{eq:psiH}, we split it into three regions as follows. We denote $\tilde g(r):=g(r)\zeta_H(r)  r^{N-1}$.
\beq
\label{eq:Rsplit}
\begin{aligned}
\psi_H(\rho_1)-\psi_H(\rho_2)&=(R1)+(R2)+(R3)\\
(R1)&=\int_{a_\it}^{\rho_1}\l(\rho_1^{2s-N}\Phi\l(\frac{r}{\rho_1}\r)-\rho_2^{2s-N}\Phi\l(\frac{r}{\rho_2}\r)\r)  \tilde g(r)\d r\\
(R2)&=\int_{\rho_1}^{\rho_2}\l(\rho_1^{2s-N}\Phi\l(\frac{r}{\rho_1}\r)-\rho_2^{2s-N}\Phi\l(\frac{r}{\rho_2}\r)\r) \tilde  g(r) \d r\\
(R3)&=\int_{\rho_2}^{b_\it}\l(\rho_1^{2s-N}\Phi\l(\frac{r}{\rho_1}\r)-\rho_2^{2s-N}\Phi\l(\frac{r}{\rho_2}\r)\r)  \tilde g(r)\d r.
\end{aligned}
\eeq
Next we will estimate $(R1)-(R3)$ separately (though the developments for $(R1)$ and $(R3)$ will be parallel).\\ 

\uline{Bound for $(R1)$.}
Define $\curly{C}_H:=[a_H,b_H]\times [a_\it,b_\it]$ and note that $\curly{C}_H\cc\curly{C}$ from \eqref{eq:Cdefn}. We use the first equation in \eqref{eq:PhiK} in \eqref{eq:PhiK} to get 
$$
\begin{aligned}
 (R1) 
 =&\int_{a_\it}^{\rho_1} \l(\rho_1^{2s-N}\vp_1(\rho_1,r)-\rho_2^{2s-N}\vp_1(\rho_2,r)\r) \tilde g(r)\d r\\
 &+\int_{a_\it}^{\rho_1} K_{2s}(\rho_1-r)\l(\rho_1^{2s-N}\vp_2(\rho_1,r)-K_{2s}(\rho_2-r)\rho_2^{2s-N}\vp_2(\rho_2,r)			\r)  \tilde g(r)\d r
 \end{aligned}
$$
Since $\curly{C}\cc \curly{C}_H$, the $\vp_i$ and all power functions $(\rho,r)\mapsto \rho^\beta r^\gam$ ($\beta,\gam\in \R$) are bounded and have bounded derivatives on $\curly{C}_H$.  In particular, they are uniformly Lipschitz continuous by the mean value theorem. Since $\tilde g$ is locally integrable, the triangle inequality gives
$$
\begin{aligned}
|(R1)|\leq C|\rho_1-\rho_2|+C'\int_{a_\it}^{\rho_1}\l|K_{2s}(\rho_1-r)-K_{2s}(\rho_2-r)\r|  |\tilde g(r)|\d r.
\end{aligned}
$$
The first term is even Lipschitz continuous. After estimating the upper bound in the integral by $\rho_1\leq b_\it$, the second term is the special case $N=1$ of what was treated in the conclusion of Theorem \ref{thm:Holder1}, see \eqref{eq:proof14}-\eqref{eq:rho}. The claimed H\"older/almost Lipschitz estimate thus holds for $(R1)$.

 
\uline{Bound for (R3).} The argument is completely analogous to that for (R1). The only difference is that one uses the second line in \eqref{eq:PhiK} instead of the first line.

\uline{Bound for (R2).} Here we need more precise information than \eqref{eq:PhiK} (which does not give any relation between the $\vp_i$). Instead we return to \eqref{eq:Phicases}, assuming $s\neq 1/2$ for now. We have
$$
\begin{aligned}
(R2)=&\int_{\rho_1}^{\rho_2}\l(\rho_1^{2s-N}\Phi\l(\frac{r}{\rho_1}\r)-\rho_2^{2s-N}\Phi\l(\frac{r}{\rho_2}\r)\r) \tilde g(r)\d r\\
=&\int_{\rho_1}^{\rho_2}
\l(\rho_1^{2s-N}\l(\frac{r}{\rho_1}\r)^{2s-N}\phi_1\l(\frac{\rho_1}{r}\r)-\rho_2^{2s-N}\phi_1\l(\frac{r}{\rho_2}\r)\r) \tilde g(r)\d r\\
&+\int_{\rho_1}^{\rho_2} 
\l(\l|\rho_1-r\r|^{2s-1}r^{1-N}\phi_2\l(\frac{\rho_1}{r}\r)-\l|\rho_2-r\r|^{2s-1}\rho_2^{1-N}\phi_2\l(\frac{r}{\rho_2}\r)\r)  \tilde g(r)\d r.
\end{aligned}
$$
We recall that the $\phi_i$ and all powers in $\rho_1,\rho_2$ and $r$ are smooth, bounded functions with uniformly bounded derivatives for the values of $\rho_1,\rho_2$ and $r$ considered here. Hence, by the mean value theorem and the key fact that $|r-\rho_i|\leq|\rho_1-\rho_2|$ ($i=1,2$) holds for all $\rho_1<r<\rho_2$, we have
$$
|(R2)|\leq C|\rho_1-\rho_2|+C'\int_{\rho_1}^{\rho_2}	\l|	\l|\rho_1-r\r|^{2s-1}-\l|\rho_2-r\r|^{2s-1}\r| |\tilde g(r)|\d r.
$$
A similar argument applies to $s=1/2$ (Proposition \ref{prop:Phi} also yields a formula similar to \eqref{eq:Phicases} for $s=1/2$). The upshot is that for all $0<s<1$,
$$
\begin{aligned}
|(R2)|\leq &C|\rho_1-\rho_2|+C'\int_{\rho_1}^{\rho_2}	\l|	K_{2s}(\rho_1-r)-K_{2s}(\rho_2-r)\r| |\tilde g(r)|\d r\\
\leq& C|\rho_1-\rho_2|+C'\int_{a_\it}^{b_\it}	\l|	K_{2s}(\rho_1-r)-K_{2s}(\rho_2-r)\r| |\tilde g(r)|\d r
\end{aligned}
$$
The argument the concludes, as for $(R1)$, by following the conclusion of Theorem \ref{thm:Holder1}.
We have thus shown that $(R1)-(R3)$ separately satisfy the claimed H\"older (respectively almost Lipschitz) continuity. By \eqref{eq:Rsplit}, this proves statements (I) and (II.A) in Theorem \ref{thm:radial}.\\

\dashuline{Proof of statement (II.B) for radial $\psi$.} Recall that it suffices to prove the claimed H\"older/almost Lipschitz continuity for $\psi_H$ given in \eqref{eq:psiH}. This requires no new ideas, so we omit the details. Basically, one first shows that \eqref{eq:psiH} is differentiable in $\rho$ by a similar dominated convergence argument as in the proof of Theorem \ref{thm:Holder1} (II.B). After computing the $\frac{\d}{\d \rho}$ derivative of \eqref{eq:psiH}, one observes that, up to smooth functions, it is given by convolution with $K_{2s-1}$. Hence, one can repeat the same arguments as in the proof of statements (I) and (II.A) above.
\e{proof}

\subsection{Removing the assumption that $\psi$ is radial}
We have already shown the claim for all $N\geq 2$ if $\psi$ is radial, or equivalently if $\psi\in\curly{H}_0$. We now prove that the case $\psi\in \curly{H}_l$ for $l>0$ can be reduced to the case $l=0$. This is based on a known (but not entirely standard) trade-off between $l$ and dimension $N$: Taking Fourier transforms, we will see that, up to power functions which are smooth away from the origin, $l$ and $N$ only enter in the form $l+\frac{N-2}{2}$, see \eqref{eq:trade} below. Thus, one can set $l=0$ at the price of increasing $N$.

\be{proof}[Proof of Theorem \ref{thm:radial} for $l>0$]
Let $\psi\in\curly{H}_l$ for some $l>0$. That is, writing $x=r\om_x$ with $r>0$ and $\om_x\in\S^{N-1}$,
\beq
\label{eq:psiexp}
\psi(x)=\Psi(r)Y_l(\om_x)
\eeq
with $\Psi\in L^2(\set{R}_+;r^{N-1}\d r)$. We recall that the Fourier transform leaves each subspace $\curly{H}_l\subset L^2(\R^N)$ invariant and reduces to the Fourier-Bessel transform on it, see Theorem 3.10 in chapter IV of \cite{SteinWeiss}.
Indeed, letting $\xi=k\om_\xi$ with $k>0$ and $\om_\xi\in\S^{N-1}$, \eqref{eq:psiexp} gives
\begin{align}
\nonumber
     \ft\psi(\xi) =& i^{-l} \curly{F}_{l,N} \Psi(k)\, Y_{l}(\om_\xi)\\
     \label{eq:FB2}
         \curly{F}_{l,N} \Psi(k) :=& k^{1-N/2} \int_0^\it
    \curly{J}_{l+\frac{N-2}{2}}(rk) r^{N/2}\Psi(r)  \d r
\end{align}
Here, $\curly{J}_{l+\frac{N-2}{2}}$ denotes the Bessel function of the first kind. Recall that $V$ is radial, so that $V\psi(x)=V(r)\Psi(r)Y_l(\om_x)$. Hence, taking the Fourier transform of
$(-\De)^s \psi+ V\psi=0$ gives
\beq
\label{eq:FBeq}
    k^{2s}\, \curly{F}_{l,N} \Psi(k) + \curly{F}_{l,N}(V \Psi)(k)=0
\eeq
for almost every $k>0$. The crucial realization is that, looking at \eqref{eq:FB2}, one has
\beq
\label{eq:trade}
\curly{F}_{l,N}\phi(k)=|k|^{l}\curly{F}_{0,N+2l}(|\cdot|^{-l}\phi),\quad \forall \phi\in L^2(\set{R}_+;r^{N-1}\d r).
\eeq
Applying this to \eqref{eq:FBeq}, we obtain
\beqs
    k^{2s} \curly{F}_{0,N+2l} (|\cdot|^{-l}\Psi)(k) + \curly{F}_{0,N+2l} (V|\cdot|^{-l}\Psi)(k)=0,
\eeqs
for all $k>0$. By unitarity of the Fourier transform on the invariant subspace $\curly{H}_0$, we have
\beqs
    (-\De)^s \tilde\psi + V\tilde\psi=0,\qquad \text{in } \R^{N+2l}.
\eeqs
where we defined the \emph{radial} function $\tilde\psi(x):=r^{-l}\Psi(r)$ on $\R^{N+2l}$. Our assumption \eqref{eq:reg} on $\psi$ is designed exactly such that $g:=-V\tilde\psi$ satisfies the same assumption but with $l=0$ and $N$ replaced by $N+2l$, i.e.\
$$
g\in 
\be{cases}
\curly{L}_{1/(2s)}(\R_+;r^{N+2l-1}\d r),\quad &\text{if } 0<s\leq 1/2\\
\jap{\cdot}^{N+2l-2s}L^1(\R_+;r^{N+2l-1}\d r),\quad &\text{if } 1/2<s<1.
\e{cases}
$$
Since $\tilde{\psi}\in\curly{H}_0$ is radial, we can apply Lemma \ref{lm:step1rad} to get $\psi\in L^\it(A_\it)$ with $A'\cc A_\it\cc A$. 

\dashuline{Proof of statements (I) and (II.A).}
The arguments in the previous sections yield that $\tilde\psi\in \tilde C^{0,\al_1}([a',b'])$, where $a'=\min\setof{|y|}{y\in A'}$ and $b'=\max\setof{|y|}{y\in A'}$. Since $a'>0$, all power functions are smooth on an open neighborhood of $[a',b']$ and we obtain $$\Psi=(\cdot)^l\tilde\psi\in \tilde C^{0,\al_1}([a',b']).$$ We take points $x_1,x_2\in A'$ and let $\al_1<1$ for definiteness. By \eqref{eq:psiexp} we have
$$
|\psi(x_1)-\psi(x_2)|\leq \|Y_l\|_{L^\it(\S^{N-1})} ||x_1|-|x_2||^{\al_1} +\|\Psi\|_{L^\it(A')} \|\nabla_\om Y_l\|_{L^\it(\S^{N-1})} \l|\frac{x_1}{|x_1|}-\frac{x_2}{|x_2|}\r|.
$$
Here we used that the spherical harmonics are polynomials (in particular they are smooth) on the compact set $\S^{N-1}$. Recall also that $\|\Psi\|_{L^\it(A')}<\it$ by Lemma \ref{lm:step1rad}. Since $|x_1|$ and $|x_2|$ are bounded away from zero, several applications of the triangle inequality give 
$$
|\psi(x_1)-\psi(x_2)|\leq C||x_1|-|x_2||^{\al_1}+ C'||x_1|-|x_2||+C''|x_1-x_2|\leq C''|x_1-x_2|^{\al_1}.
$$
In the second step we used that we may restrict to $|x_1-x_2|<1$ at the price of increasing the constant. The same argument applies when $\al_1=1$ and we have an additional logarithm.

\dashuline{Proof of statement (II.B).}
By the results of the previous sections, we get $\Psi\in \tilde C^{1,\al_1-1}([a',b'])$. We compute the gradient of $\psi$ given by \eqref{eq:psiexp} in spherical coordinates and get, for $x=r\om\in A'$,
$$
\nabla\psi(x)=\Psi'(r)Y_l(\om) \om + \frac{\Psi(r)}{r} \l(\nabla_\om Y_l(\om)-\om(\om\cdot\nabla_\om)Y_l(\om)\r).
$$
Since $Y_l$ is bounded and $\Psi\in \tilde C^{1,\al_1-1}([a',b'])$, the first term is indeed in $\tilde C^{0,\al_1-1}(A')$. The second term is in fact differentiable in $r$ and $\om$  on $A'$. (We remark that the term $\nabla Y_l(\om)-\om(\om\cdot\nabla)Y_l(\om)$ is known as a ``vector spherical harmonic''.) Hence, letting $\al_1-1<1$ for definiteness, we have
$$
|\nabla\psi(x_1)-\nabla\psi(x_2)|\leq C ||x_1|-|x_2||^{\al_1-1}+ C'||x_1|-|x_2||+C'' \l|\frac{x_1}{|x_1|}-\frac{x_2}{|x_2|}\r|
\leq |x_1-x_2|^{\al_1-1}.
$$
This finishes the proof of Theorem \ref{thm:radial}.
\e{proof}

\section{Proof of Theorem \ref{thm:Holder1local}}

\subsection{Non local Leibniz rule}

\be{prop}[Non local Leibniz rule]
\label{prop:nlleibniz}
Let $\zeta,\vp\in C_0^\it(\set{R}^N)$.
\be{enumerate}[label=(\roman*)]
	\item If $0<s< 1/2$, then for all $x\in\R^N$,
	\beq
		 (-\De)^s(\zeta\vp)(x)= \zeta(x)(-\De)^s \vp(x) - \int \frac{\zeta(x)-\zeta(y)}{|x-y|^{N+2s}} \vp(y) \d y.
	\eeq
\item If $1/2\leq s< 1$, then for all $x\in\R^N$,
	\beq
	\label{eq:nlleibniz2}
	\begin{aligned}
		 &(-\De)^s(\zeta\vp)(x)\\
		 &= \zeta(x)(-\De)^s \vp(x) + \l((-\De)^s \zeta\r)(x) \vp(x) - \int \frac{(\zeta(x)-\zeta(y))(\vp(x)-\vp(y))}{|x-y|^{N+2s}} \d y.
\end{aligned}	
	\eeq
\e{enumerate}
\e{prop}

\be{proof}
This is obtained by simple algebra from the integral formula \eqref{eq:FLdefn} for $(-\De)^s$ acting on Schwartz functions. 
\e{proof}

We can use Proposition \ref{prop:nlleibniz} to extend the equation from $\Om$ to the whole space, up to a localization error.

\be{cor}
\label{cor:FSEwhole}
We make the same assumptions as in Theorem \ref{thm:Holder1local}, in particular
\beq
\label{eq:FSEOmega}
    (-\De)^s\psi + V\psi=f\quad \text{on }  \Om
\eeq
holds in the distributional sense. Let $\zeta\in C_0^\it(\R^N)$ be a cutoff function with $\supp\, \zeta\cc \Om$ and $B\cc \{\zeta=1\}$. Then:

\be{enumerate}[label=(\roman*)]
	\item If $0<s< 1/2$, we have
	\beq
\label{eq:FSEwhole}
    (-\De)^s(\zeta \psi) = \zeta(-V\psi+f) + \curly{E}_1 \quad \text{on }  \set{R}^N,
\eeq
	in the distributional sense. Here, the localization error $\curly{E}_1$ is defined by
	\beq
	\label{eq:E1defn}
	 \curly{E}_1(x)= -\int \frac{\zeta(x)-\zeta(y)}{|x-y|^{N+2s}} \psi(y) \d y.
	\eeq
\item If $1/2\leq s<1$, we have
	\beq
	\label{eq:FSEwhole2}
		 (-\De)^s(\zeta\vp)= \zeta(-V\psi+f)+((-\De)^s\zeta) \psi + \curly{E}_2		 
	\eeq
		in the distributional sense. Here, the localization error $\curly{E}_2$ is defined by
	\beq
	\label{eq:E2defn}
	 \curly{E}_2(x)= -\int \frac{(\zeta(x)-\zeta(y))(\psi(x)-\psi(y))}{|x-y|^{N+2s}}\d y
	 \eeq
\e{enumerate}
\e{cor}

\be{proof}
Since $\psi$ is a distributional solution to \eqref{eq:FSEOmega} and $\supp\,\zeta\cc \Om$, we get
$\zeta\psi\in \jap{\cdot}^{N+2s}L^1(\set{R}^N)$
 and $\zeta f,\zeta V\psi\in L^1_{\mathrm{loc}}(\set{R}^N)$. We will prove momentarily, in Lemma \ref{lm:locerror}, that
 $\curly{E}_1, \curly{E}_2\in L^1(\set{R}^N)$. Thus, it remains to check \eqref{eq:FSEdistribution}. 
 
 Let $\vp\in C_0^\it(\set{R}^N)$ and $0<s<1/2$. By Proposition \ref{prop:nlleibniz} and Fubini's theorem, we have
\beq
\label{eq:key}
\scp{ (-\Del)^s\vp}{\zeta \psi} = \scp{\zeta  (-\Del)^s\vp}{\psi}=\scp{(-\Del)^s(\zeta\vp)}{\psi} - \scp{\vp}{\curly{E}_1}.
\eeq
 We observe that $\zeta\vp\in C_0^\it(\Om)$ and use \eqref{eq:FSEOmega} to conclude
\beqs
   \scp{(-\Del)^s(\zeta\vp)}{\psi} = \scp{\vp}{\zeta(-V\psi+f)}.
\eeqs
Since $\vp\in C_0^\it(\set{R}^N)$ was arbitrary, this proves (i). For (ii), the claim follows by an analogous argument.
\e{proof}


 \subsection{Controlling the localization error}
 
 The non local Leibniz rule allows us to extend the equation from $\Om$ to an equation on $\R^N$ with new $\tilde V$ and $\tilde f$ (e.g.\ \eqref{eq:FSEwhole} has $\tilde f = \zeta f+\curly{E}_1$ and $\tilde V = \zeta V$). To apply Theorem \ref{thm:Holder1}, we need the appropriate integrability of the localization errors $\curly{E}_1, \curly{E}_2$.

The following lemma says that $\curly{E}_1$ is bounded on $B$ and is therefore unproblematic. It also says that $\curly{E}_2$ has the same local integrability as $\psi$, which will improve locally as we perform the bootstrap procedure, culminating in a locally bounded $\curly{E}_2$ as well. 
 
\be{lm}
\label{lm:locerror}
Let $\psi, \zeta$ be as in Corollary \eqref{cor:FSEwhole}.
\be{enumerate}[label=(\roman*)]
	\item If $0<s< 1/2$,  then $\curly{E}_1\in L^\it(B)\cap L^1(\set{R}^N)\cap \curly{L}_{N/(2s)}(\R^N)$.
	\item If $1/2\leq s<1$ and $\psi\in L^q(\tilde B)$ for some $\tilde B\cc B$, then $\curly{E}_2 \in L^q(\tilde B)\cap L^1(\set{R}^N)\cap \curly{L}_{N/(2s)}(\R^N)$.
\e{enumerate}	
\e{lm} 
 
 
 \be{proof}
 \dashuline{Proof of statement (i).} Consider \eqref{eq:E1defn}. For $x\in B$, we have $\zeta(x)=1$ and so $\zeta(x)-\zeta(y)=0$ unless $y\notin \{\zeta=1\}$. Since $B\cc \{\zeta=1\}$, this implies that $|x-y|\geq C$ holds whenever $\zeta(x)-\zeta(y)\neq 0$. Together with $|x|\leq C'$, this implies $|x-y|\geq C''\jap{y}$ and so
 \beqs
  |\curly{E}_1(x)|\leq C \int_{\{\zeta\neq1\}} \frac{|\psi(y)|}{\jap{y}^{N+2s}}\d y= C'
 \eeqs
 holds for all $x\in B$. The second is finite by Definition \ref{defn:distributional} of a distributional solution. We conclude that $\curly{E}_1\in L^\it(B)$. Now let $B',B''$ be balls such that $B\cc \supp\,\zeta\cc B'\cc B''$ and let $x\in B'$. We apply the mean value theorem to $\zeta(x)-\zeta(y)$ to find
 \beqs
  |\curly{E}_1(x)|\leq C \int_{B''} \frac{|\psi(y)|}{|x-y|^{N+2s-1}}\d y 
  + C' \int_{(B'')^c} \frac{|\psi(y)|}{\jap{y}^{N+2s}}\d y.
 \eeqs 
 By Young's inequality, the first term lies in $L^{1+\de}(B')$ for sufficiently small $\de>0$ (as a function of $x$). The second term evaluates to a finite constant and thus lies in $L^\it(B')$. Finally, let $x\in (B')^c$. Then, $\zeta(x)=0$ and so $\zeta(x)-\zeta(y)=0$ unless $y\in \supp\, \zeta$. Since $\supp\, \zeta\cc B'$, this implies that $|x-y|\geq C\jap{x}$ and so
  \beq
  \label{eq:E1decay}
  |\curly{E}_1(x)|\leq  \frac{C}{\jap{x}^{N+2s}} \int_{\supp\,\zeta} |\psi(y)| \d y= \frac{C'}{\jap{x}^{N+2s}}
 \eeq
 for all $x\in (B')^c$. We have shown that writing $\curly{E}_1=\chi_{B'}\curly{E}_1+\chi_{(B')^c}\curly{E}_1$ the first function is in $L^{1+\de}(\R^N)$ and the second function satisfies $\chi_{(B')^c}|\curly{E}_1|\leq C\jap{x}^{-N-2s} $. Recalling Definition \ref{defn:weightedL1} of $\curly{L}_{N/(2s)}(\R^N)$, the claim follows.\\
 
 
 \dashuline{Proof of statement (ii).}  Consider \eqref{eq:E2defn}. For $x\in \tilde B$, we have $\zeta(x)=1$ and we see that $\zeta(x)-\zeta(y)\neq 0$ gives a lower bound on $|x-y|$ and so
  \beqs
  |\curly{E}_2(x)|\leq C|\psi(x)| + C'.
 \eeqs
 for all $x\in \tilde B$. So clearly $\psi \in L^q(\tilde B)$ implies $\curly{E}_2 \in L^q(\tilde B)$. Let $B',B''$ be balls such that $B\cc \supp\,\zeta\cc B'\cc B''\cc \Om$. For $x\in (B')^c$, the same argument that gave \eqref{eq:E1decay} yields
   \beq
  \label{eq:E2decay}
  |\curly{E}_2(x)|\leq  C \frac{|\psi(x)|+1}{\jap{x}^{N+2s}}
 \eeq
and this is integrable since $\psi$ is a distributional solution. Considering Definition \ref{defn:weightedL1} of $\curly{L}_{N/(2s)}(\R^N)$, it remains to show that $\curly{E}_2\in L^{1+\de}(B')$ for some $\de>0$. For $x\in B'$, we have
$$
\begin{aligned}
|\curly{E}_2(x)|&\leq \int_{B''}\l|\frac{(\zeta(x)-\zeta(y))(\psi(x)-\psi(y))}{|x-y|^{N+2s}}\r|\d y
+\int_{(B'')^c}\frac{|\psi(x)|+|\psi(y)|}{\jap{y}^{N+2s}}\d y\\
&\leq \int_{B''}\l|\frac{\psi(x)-\psi(y)}{|x-y|^{N+2s-1}}\r|\d y
+C|\psi(x)|+C'.
\end{aligned}
$$
Recall our assumption that $\psi\in W^{s',1+\de}(\Om)$ for some $s'>2s-1$ and $\de>0$. In particular, $\psi\in L^{1+\de}(\Om)$ and so it suffices to consider the integral in the last expression (as a function of $x$). By Jensen's inequality and $B'\cc B''\cc \Om$,
$$
\int_{B'}\l(\int_{B''}\l|\frac{\psi(x)-\psi(y)}{|x-y|^{N+2s-1}}\r|\d y\r)^{1+\de}\d x
\leq C\int_\Om\int_\Om \l|\frac{\psi(x)-\psi(y)}{|x-y|^{N+2s-1}}\r|^{1+\de}\d x \d y
\leq C\|\psi\|_{H^{\tilde s,1+\de}}
$$
with $\tilde s=2s-1+c\de$ for some $c>0$ independent of $\de$. By decreasing $\de$ if necessary (the fractional Sobolev spaces on $\Om$ are nested sets), we can set $\tilde s =s'$ and we are done.
 \e{proof}
 



\subsection{Conclusion}
We have extended the equation to $\set{R}^N$ and have control over the localization error. We are now ready to give the

\be{proof}[Proof of Theorem \ref{thm:Holder1local}]
\dashuline{Proof of statement (I).} Let $0<s<1/2$. By \eqref{eq:FSEwhole}, $\zeta\psi$ satisfies  
\beq
\label{eq:FSEmod}
(-\De)^s(\zeta\psi)=-\tilde V\psi+\tilde f
\eeq
 with $\tilde V=\zeta V$ and $\tilde f=\zeta f+\curly{E}_1$. Since $\zeta$ is smooth and $\curly{E}_1\in L^\it(B)\cap \curly{L}{N/(2s)}(\R^N)$ by Lemma \ref{lm:locerror}, all the assumptions in Theorem \ref{thm:Holder1} are met and we obtain the claimed H\"older regularity for $\zeta\psi$ on $B'$ and hence also for $\psi$ (as $\zeta= 1$ on $B$). 
 
 

\dashuline{Proof of statement (II).} Let $1/2\leq s\leq 1$. By \eqref{eq:FSEwhole2}, $\zeta\psi$ satisfies \eqref{eq:FSEmod}, but now with 
\beqs
\tilde V=\zeta V+(-\De)^s\zeta,\qquad \tilde f=\zeta f+\curly{E}_2.
\eeqs
We have $\tilde V,\zeta f\in L^p(B)$ since $\zeta$ is smooth and $\curly{E}_2\in \curly{L}_{N/(2s)}(\R^N)$ by Lemma \ref{lm:locerror}.
However, we no longer have that the localization error $\curly{E}_2$ is bounded on $B$. Instead, Lemma \ref{lm:locerror} implies that $\curly{E}_2$ inherits the local integrability of $\psi$.

This allows us to run a modified bootstrap procedure: We \emph{interlace} each iteration step of the previous bootstrap argument, see the proof of Lemma \ref{lm:step1}, with an application of Lemma \ref{lm:locerror} to get that the local integrability of the localization error $\curly{E}_2$ has also improved. The details are left to the reader. This gives boundedness of $\psi$ on a smaller ball $B_\it\cc B$ and the same argument as in Section \ref{ssect:conclusion} yields statement (II).
\e{proof}


\section*{Acknowledgments}
It is a pleasure to thank Rupert Frank for suggesting the problem and for several helpful discussions. We thank Tianling Jin for useful remarks on a draft version of this paper.

\be{appendix}
\section{Proof of Proposition \ref{prop:finite}}
\label{app:propfinite}
Since $\psi\in L_{\mathrm{loc}}^1(\Om)$, it suffices to prove that $\scp{(-\Del)^s\vp}{\psi}$ makes sense when $s<1$. The idea is to use \eqref{eq:FLdefn} for large $x$ to get the explicit polynomial decay rate of $(-\De)^s\vp(x)$. Let $B$ be an open ball such that $\supp\,\vp\cc B$. By \eqref{eq:FLdefn}, we have
\beq
\label{eq:termatinfty}
		\scp{(-\Del)^s\vp}{\psi}	= \int_{B} \ol{(-\Del)^s \vp(x)} \psi(x) \dx		
	 - \int_{B^c}\l(\int_{\R^N} \frac{\vp(y)}{|x-y|^{N+2s}} \d y\r)\psi(x) \d x.
\eeq
Here we could drop the principal value from \eqref{eq:FLdefn} because $\textnormal{dist}(B^c,\supp\,\vp)>0$. It suffices to prove that
\beq
\label{eq:weightedcontinuity}
    \scp{(-\Del)^s\vp}{\psi} \leq C_\vp \|\jap{\cdot}^{-N-2s}\psi\|_{L^1(\R^N)}<\it.
\eeq
For the first term in \eqref{eq:termatinfty}, finiteness follows directly from the observations that $\psi\in L^1_{\mathrm{loc}}(\R^N)$ and that $(-\De)^s\vp$ is a smooth function. (Smoothness can be seen from the Sobolev embedding theorem and the fact that the Fourier transform of $(-\De)^s\vp$ decays faster than any polynomial.) For the second term in \eqref{eq:termatinfty}, we make the general observation that $|y|\leq C$ and $|x-y|\geq C'$ implies that $|x-y|\geq C''\jap{x}$ with $C''$ depending only on $C'$ and $C''$. Then \eqref{eq:weightedcontinuity} follows.
\qed

\section{Proof of Proposition \ref{prop:riesz}}
\label{app:propriesz}
We separate the proof into two Lemmas. Together they imply Proposition \ref{prop:riesz} by linearity. Recall that $0<s\leq \min\{1,N/2\}$.

\be{lm}
\label{lm:riesz}
Let $f\in \curly{L}_{N/(2s)}(\R^N)$. Then, in the distributional sense of Definition \ref{defn:distributional},
\beqs
(-\De)^s(k_{2s} * f)=f \quad \text{on } \set{R}^N.
\eeqs
\e{lm}

\be{lm}
\label{lm:uniqueness}
Suppose that $(-\De)^s \psi_0=0$ holds in the distributional  sense of Definition \ref{defn:distributional}. Then:
 \be{enumerate}[label=(\roman*)]
\item If $s\leq1/2$, there exists $w\in\set{C}$ such that $\psi_0(x)= w$ for a.e.\ $x\in\set{R}^N$.
\item If $1/2< s\leq 1$, there exist $w,w_1,\ldots,w_N\in\set{C}$ such that $\psi_0(x) = w+ \sum_{j=1}^N w_j x_j$
for a.e.\ $x\in\set{R}^N$.
\e{enumerate}
\e{lm}

\be{proof}[Proof of Lemma \ref{lm:riesz}]
Let $\vp\in C_0^\it(\set{R}^N)$ be arbitrary. Suppose for the moment that $f$ is a Schwartz function. When $N=2s=2$, the result is well-known. For $2s<N$, we will use Proposition \ref{prop:rieszkernel} and for $N=2s=1$, we will use the fact that
\beq
\mathrm{P.V.} \frac{1}{|\xi|}=-\sqrt{8\pi}(\gam_{EM}+\log|\cdot|)^\wedge (\xi)\equiv \ft{k_{1}}(\xi),
\eeq
see formula (32) on p.\ 132 of \cite{Vladimirov} (though a different convention for the Fourier transform is used there). Then, Fubini's theorem and the convolution theorem for tempered distributions, see e.g.\ Theorem IX.4 in \cite{ReedSimonII}, gives
\beq
\label{eq:schwartzproof}
    \scp{(-\De)^s\vp}{k_{2s}*f}=\scp{k_{2s}*(-\De)^s\vp}{f} =\scp{(2\pi)^{N/2}(\ft{k_{2s}} |\cdot|^{2s}\ft{\vp})^\vee}{f} =\scp{\vp}{f}.
\eeq
Hence the claim holds if $f$ is a Schwartz function. Now let $f\in \curly{L}_{N/(2s)}(\R^N)$. Consider first the case $2s<N$ and write $L^q\equiv L^q(\R^N)$. By the Definition \ref{defn:weightedL1} of $f\in\curly{L}_{N/(2s)}(\R^N)$, there exists $\de>0$ such that $f=f_1+f_2$ with $f_1\in L^{1+\de}$ and $f_2\in L^{N/(2s)-\de}$. We claim that
\beq
\label{eq:approxremains}
\begin{aligned}
|\scp{(-\De)^s\vp}{k_{2s}*f_1}|\leq C_\vp \|f_1\|_{L^{1+\de}},\quad |\scp{(-\De)^s\vp}{k_{2s}*f_2}|\leq& C_\vp \|f_2\|_{L^{N/(2s)-\de}}.
\end{aligned}
\eeq
This will imply the lemma after approximating $f_1,f_2$ with Schwartz functions. To prove \eqref{eq:approxremains}, we note that  $|(-\De)^s\vp(x)|\leq C_\vp\jap{x}^{-N-2s}$ for any $\vp\in C_0^\it(\R^N)$. Indeed, this is trivial for $s=1$ and for $s<1$ it follows by combining the fact that $(-\De)^s\vp$ is smooth and its decay at infinity being $\jap{x}^{-N-2s}$ by \eqref{eq:FLdefn}. Then \eqref{eq:approxremains} follows from  the Hardy-Littlewood-Sobolev inequality, see e.g.\ \cite{LiebLoss}. 
 
Finally, let $2s=N$. Our assumption $f\in \curly{L}_1(\R^N)$ means that $f\in L^{1+\de}\cap \jap{\cdot}^{-\de}L^1$ for some $\de>0$. 
Similarly as \eqref{eq:approxremains} above, we show that
\beq
\label{eq:approxremains2}
|\scp{(-\De)^s\vp}{k_{2s}*f_1}|\leq C_\vp \max\{\|f\|_{L^{1+\de}},\|\jap{\cdot}^{\de}f\|_{L^1}\}.
\eeq
Indeed, for every $\eps>0$ we have
$$
\begin{aligned}
|k_{2s}(x-y)|\leq &C |\log|x-y||+C'\leq C'' \max\{|x-y|^{-\eps},|x-y|^\eps\}\\
\leq& C''' \max\{|x-y|^{-\eps},\jap{x}^\eps+\jap{y}^\eps\}.
\end{aligned}
$$ 
 Choosing $\eps=\min\{s,\de\}$, yields \eqref{eq:approxremains2} by the Hardy-Littlewood-Sobolev inequality. The approximation argument to go from \eqref{eq:approxremains2} to the main claim for $2s=N$ is less trivial than in the case $2s<N$ (one has to produce a sequence of Schwartz functions converging to $f$ in \emph{both} $L^{1+\de}$ and $\jap{\cdot}^{-\de}L^1$). Apart from the usual trick of restricting to bounded and compactly supported $f$, the key idea is to approximate $\jap{\cdot}^\de f\in L^1$ by the usual sequence of mollified $\phi_n\in C_0^\it(\R^N)$ and to observe that $\tilde \phi_n:=\jap{\cdot}^{-\de}\phi_n$ are still Schwartz functions. The standard arguments give $\tilde \phi_n\to f$ in $L^{1+\de}(\R^N)$ and this proves the lemma.
\e{proof}

\be{proof}[Proof of Lemma \ref{lm:uniqueness}]
\dashuline{Step 1.} We first prove that the distributional support of $\ft{\psi_0}$ is $\{0\}$. To this end, take any $\eta\in C_0^\it(\set{R}^N\setminus\{0\})$. By the definition of the Fourier transform on tempered distributions, we have
\beqs
    \scp{\eta}{\ft{\psi_0}}  =\scp{|\cdot|^{2s}|\cdot|^{-2s}\eta}{\ft{\psi_0}} =\scp{(-\De^s)\l(|\cdot|^{-2s}\eta\r)^\vee}{\psi_0}.
\eeqs
If $\l(|\cdot|^{-2s}\eta\r)^\vee$ were a $C_0^\it$-function, the last expression would vanish by our assumption that $(-\De)^s\psi_0=0$. Note that
 $|\cdot|^{-2s}\eta\in C_0^\it(\set{R}^N\setminus\{0\})$ since $\eta$ vanishes near the origin where $|\cdot|^{-2s}$ is non smooth. However, this
only tells us that its inverse Fourier transform is a Schwartz function. (In fact, by a version of the uncertainty principle, it cannot be
compactly supported.)

We proceed by an approximation argument. By denseness, we can find $(\vp_n)_{n\geq1}\subset
C_0^\it(\set{R}^N)$ such that $\vp_n\rightarrow \l(|\cdot|^{-2s}\eta\r)^\vee$ in the Schwartz space topology.
Since the Fourier transform is an isometry on Schwartz space, we also have $\ft{\vp}_n\rightarrow |\cdot|^{-2s}\eta$
 in the topology of $\curly{S}(\set{R}^N)$. Note however that we do not have $|\cdot|^{2s}\ft{\vp}_n\rightarrow \eta$ in the Schwartz space topology, since $|\cdot|^{2s}$ is non-smooth at the origin and the approximating sequence $\ft{\vp}_n$ need not vanish in an open neighborhood of the origin just because its limit $\eta$ does.
 
 Nonetheless, we can modify the approximating sequence to improve the convergence by hand. For this, we acknowledge a helpful discussion on mathoverflow.net \cite{MO}. For $k\geq 0$, let
  \beqs
    \eps_{n,k}: = \ft{\vp}^{(k)}_n(0)=(2\pi)^{-3/2} \int (ix)^k \vp_n(x) \dx,
  \eeqs
 which converges to $\eta^{(k)}(0)=0$ as $n\rightarrow \it$. Let $\zeta_k\in C_0^\it(\set{R}^N)$ be such that
 \beqs
    (2\pi)^{-3/2}\int (ix)^k\zeta_k(x)\d x =1.
 \eeqs
 For any integer $M>0$, we define the modified approximating sequence
 \beqs
    \tilde{\vp}_n:= \vp_n - \sum_{k=0}^{M+1} \eps_{n,k} \zeta_k,
 \eeqs
 which is clearly in $C_0^\it(\set{R}^N)$. Moreover, $\ft{\tilde{\vp}}_n^{(k)}(0)=0$ for all $0\leq k\leq M+1$ and consequently
 $|\cdot|^{2s}\ft{\tilde{\vp}}_n\rightarrow \eta$ in $C^{M}$.

 By the Hausdorff-Young inequality, $(-\De)^s\tilde \vp_n\rightarrow \eta$ in the weighted
  space $\jap{\cdot}^{-M} L^\it(\set{R}^N)$. Recall that by Definition \ref{defn:distributional} of a distributional solution,
  $\psi_0\in \jap{\cdot}^{N+2s}L^1(\set{R}^N)$. Thus, if we choose $M>N+2s$, we have the right duality product to conclude
  \beqs
    \scp{\eta}{\ft{\psi}_0}=\scp{\check\eta}{\psi_0} =\lim_{n\rightarrow \it} \scp{(-\De)^s\tilde{\vp}_n}{\psi_0}=0.
  \eeqs
Since $\eta$ was arbitrary, the distributional support of $\ft{\psi}_0$ is indeed $\{0\}$.\\

\dashuline{ Step 2.} Since the distributional support of $\ft{\psi}_0$ is $\{0\}$, Theorems 6.24 and 6.25 in \cite{RudinFA} imply that there exists a non-negative integer $K$, such that
\beq
\label{eq:psi0sum}
    \ft{\psi}_0 = \sum_{|\al|\leq K} c_\al D^\al \delta.
\eeq
Here, the sum runs over multi-indices $\al$ and $c_\al$ are constants. Note that $
    (D^\al\de)^\vee(x) = (2\pi)^{-3/2} (ix)^\al$ holds
 in the sense of tempered distributions. Recall that distributional solutions need to satisfy $\psi\in \jap{\cdot}^{N+2s}L^1(\set{R}^N)$ by
 Definition \ref{defn:distributional}. Therefore, only terms with $|\al|< 2s$ contribute to the sum \eqref{eq:psi0sum}. This proves Lemma \ref{lm:uniqueness} and hence Proposition \ref{prop:riesz}.
\e{proof}

\section{Proof of Proposition \ref{thm:optimality}}
\label{app:optimality}

The key observation is

\be{prop}
For $0<\beta < 2s$, there exists $C(\beta)\in\set{R}\setminus\{0\}$ such that
\label{prop:scaling}
\beq
\label{eq:FLscaling}
    (-\Del)^s |x|^{\beta} = C(\beta) |x|^{\beta-2s}
\eeq
as tempered distributions. 
\e{prop}

\be{proof}
One can either use \eqref{eq:FLdefn} or, as we do here, analyticity. Let $0<z<N$, then Proposition \ref{prop:rieszkernel} yields the following equality, in the sense of tempered distributions
\beq
\label{eq:AC}
\ft{((-\De)^s k_z)}(\xi)=|\xi|^{2s-z}.
\eeq
Recalling the definition \eqref{eq:rieszdefn} of $k_z$, we see that both sides are analytic in $z$ and make sense on the larger set 
$$
\setof{z\in\C}{0<\mathrm{Re}(z)<N+2s}\setminus\{0\}.
$$
Therefore \eqref{eq:AC} extends to this set. Now \eqref{eq:FLscaling} follows by setting $z=N+\beta$ and applying the inverse Fourier transform to both sides of \eqref{eq:AC} in the sense of tempered distributions.
\e{proof}

We can now give

\be{proof}[Proof of Proposition \ref{thm:optimality}]
We shall only consider the cases (I) and (II.A) in Theorem \ref{thm:Holder1}, so $s,N,p$ are such that $\al=2s-N/p<1$. We let $V\equiv 0$.
By the translation invariance of $(-\De)^s$, we may assume that $B'$ contains the origin. Let $\zeta\in C_0^\it(\R^N)$ be a cutoff function with $B\cc \{\zeta=1\}$. Let $\eps>0$ be small enough such that $\al+2\eps <1$.
Define
\beqs
    \psi(x) := \zeta(x) |x|^{\al+\eps}
\eeqs
First, observe that $\psi$ is \emph{not} H\"older continuous of order $\al+2\eps$ at the origin. It is however bounded and compactly supported. In particular, $\psi \in \jap{\cdot}^{N+2s}L^1(\set{R}^N)$ and therefore it qualifies as a distributional solution, see Definition \ref{defn:distributional}. Trivially, it solves the fractional Schr\"odinger equation with $V\equiv 0$ and $f=(-\De)^s\psi$. The claimed optimality will follow (since $\eps>0$ is arbitrarily small) once we show that\beq
\label{eq:fremains}
f=(-\De)^s(\zeta|\cdot|^{\al+\eps})\in L^p(B)\cap \curly{L}_{N/(2s)}(\R^N).
\eeq

Let $0<s<1/2$. We use the non local Leibniz rule in the same way as in \eqref{eq:key} and then we apply Proposition \ref{prop:scaling} with $\beta=2s-N/p+\eps$ (which lies in $(0,2s)$ for small enough $\eps>0$) to get
$$
f(x)=\zeta(x)|x|^{-N/p+\eps}-\int \frac{\zeta(x)-\zeta(y)}{|x-y|^{N+2s}}|y|^{\al+\eps}\dy.
$$
Note that the integral is just the localization error $\curly{E}_1$ from \eqref{eq:E1defn} but with $\psi(y)$ replaced by $|y|^{\al+\eps}$ (which we note also lies in $\jap{\cdot}^{N+2s}L^1(\R^N)$ for sufficiently small $\eps>0$). Hence, the proof of Lemma \ref{lm:locerror} gives $\curly{E}_1\in L^\it(B)\cap L^1(\R^N)\cap \curly{L}_{N/(2s)}(\R^N)$. Since $\zeta(x)|x|^{-N/p+\eps}\in L^p_{\mathrm{loc}}(\R^N)$ and it is compactly supported, \eqref{eq:fremains} follows. In the case $1/2\leq s\leq 1$, one can similarly combine the arguments from Section 6 and Proposition \ref{prop:scaling} to get \eqref{eq:fremains}. The details and the very similar arguments in the other cases (Theorem \ref{thm:Holder1}(II) and Theorem \ref{thm:Holder1local}) are left to the reader. 
\e{proof}

\e{appendix}

\bibliographystyle{amsplain}
\bibliography{PDE}

\providecommand{\bysame}{\leavevmode\hbox to3em{\hrulefill}\thinspace}
\providecommand{\MR}{\relax\ifhmode\unskip\space\fi MR }
\providecommand{\MRhref}[2]{%
  \href{http://www.ams.org/mathscinet-getitem?mr=#1}{#2}
}
\providecommand{\href}[2]{#2}
\begin{thebibliography}{10}

\bibitem{AbramowitzStegun}
M.~Abramowitz and I.A. Stegun (eds.), \emph{{Handbook of Mathematical Functions
  with Formulas, Graphs, and Mathematical Tables}}, Dover, 1972.

\bibitem{BassKassmann05}
R.F. Bass and M.~Kassmann, \emph{{H\"older continuity of harmonic functions
  with respect to operators of variable orders}}, Comm. Partial Diff. Eq.
  \textbf{30} (2005), 1249--1259.

\bibitem{BassLevin02}
Richard~F. Bass and David~A. Levin, \emph{Harnack inequalities for jump
  processes}, Potential Anal. \textbf{17} (2002), no.~4, 375--388.

\bibitem{BCF12}
C.~Bjorland, L.~Caffarelli, and A.~Figalli, \emph{Non-local gradient dependent
  operators}, Adv. Math. \textbf{230} (2012), no.~4-6, 1859--1894.

\bibitem{CS07}
Luis Caffarelli and Luis Silvestre, \emph{An extension problem related to the
  fractional {L}aplacian}, Comm. Partial Differential Equations \textbf{32}
  (2007), no.~7-9, 1245--1260.

\bibitem{CS09}
\bysame, \emph{Regularity theory for fully nonlinear integro-differential
  equations}, Comm. Pure Appl. Math. \textbf{62} (2009), no.~5, 597--638.

\bibitem{CS14}
\bysame, \emph{The {E}vans-{K}rylov theorem for nonlocal fully nonlinear
  equations}, Ann. of Math. (2) \textbf{174} (2011), no.~2, 1163--1187.

\bibitem{CaffarelliStinga14}
Luis~A. Caffarelli and Pablo~Raúl Stinga, \emph{Fractional elliptic equations,
  caccioppoli estimates and regularity}, Annales de l'Institut Henri Poincare
  (C) Non Linear Analysis (2015), --.

\bibitem{CMS}
Ren{\'e} Carmona, Wen~Chen Masters, and Barry Simon, \emph{Relativistic
  {S}chr\"odinger operators: asymptotic behavior of the eigenfunctions}, J.
  Funct. Anal. \textbf{91} (1990), no.~1, 117--142. \MR{1054115}

\bibitem{DallAcquaFournaisSorensenStockmeyer}
A.~Dall'Acqua, S.~Fournais, T.~{\O.} S{\o}rensen, and E.~Stockmeyer,
  \emph{{Real Analyticity of Solutions to Schr\"odinger Equations Involving a
  Fractional Laplacian and other Fourier Multipliers}}, Proceedings of the
  XVIIth International Congress on Mathematical Physics, World Scientific,
  2014.

\bibitem{DallAcquaFournaisSorensenStockmeyer0}
Anna Dall'Acqua, S{\o}ren Fournais, Thomas {\O}stergaard~S{\o}rensen, and
  Edgardo Stockmeyer, \emph{Real analyticity away from the nucleus of
  pseudorelativistic {H}artree-{F}ock orbitals}, Anal. PDE \textbf{5} (2012),
  no.~3, 657--691.

\bibitem{deGiorgi}
E.~De~Giogi, \emph{{Sulla differenziabilità e l'analiticità delle estremali
  degli integrali multipli regolari}}, Memorie della Accademia delle Scienze di
  Torino. Classe di Scienze Fisiche, Matematicahe e Naturali \textbf{3} (1957),
  25--43.

\bibitem{Hitchhiker}
Eleonora Di~Nezza, Giampiero Palatucci, and Enrico Valdinoci,
  \emph{Hitchhiker's guide to the fractional {S}obolev spaces}, Bull. Sci.
  Math. \textbf{136} (2012), no.~5, 521--573.

\bibitem{FelmerQuaasTan12}
P.~Felmer, A.~Quaas, and J.~Tan, \emph{{Positive solutions of the nonlinear
  Schr\"odinger equation with the fractional Laplacian}}, {Proceedings of the
  Royal Society of Edinburgh, Section: A Mathematics} \textbf{142} (2012),
  1237--1262.

\bibitem{FrankSeiringer08gs}
R.L. Frank and R.~Seiringer, \emph{{Non-linear ground state representations and
  sharp Hardy inequalities }}, Journal of Functional Analysis \textbf{255}
  (2008), no.~12, 3407 -- 3430.

\bibitem{GargSpector14}
R.~Garg and D.~Spector, \emph{{On the role of Riesz potentials in Poisson's
  equation and Sobolev embeddings}}, arXiv:1404.1563, to appear in Indiana
  Univ. Math. J.

\bibitem{GargSpector15}
Rahul Garg and Daniel Spector, \emph{On the regularity of solutions to
  {P}oisson's equation}, C. R. Math. Acad. Sci. Paris \textbf{353} (2015),
  no.~9, 819--823.

\bibitem{GilbargTrudinger}
D.~Gilbarg and N.S. Trudinger, \emph{{Elliptic Partial Differential Equations
  of Second Order}}, 2 ed., Classics in Math., Springer, Berlin, 2001.

\bibitem{GradshteynRyzhik}
I.S. Gradstein and I.M. Ryzhik, \emph{{Table of Integrals, Series and
  Products}}, 7 ed., Elsevier/Academic Press, 2007.

\bibitem{Herbst}
Ira~W. Herbst, \emph{Spectral theory of the operator
  {$(p^{2}+m^{2})^{1/2}-Ze^{2}/r$}}, Comm. Math. Phys. \textbf{53} (1977),
  no.~3, 285--294.

\bibitem{JLX}
T.~Jin, Y.~Li, and J.~Xiong, \emph{{The Nirenberg problem and its
  generalizations: A unified approach}}, arXiv:1411.5743.

\bibitem{JS15}
T.~Jin and L.~Silvestre, \emph{{Hölder gradient estimates for parabolic
  homogeneous p-Laplacian equations}}, arXiv:1505.05525.

\bibitem{KRS14}
Moritz Kassmann, Marcus Rang, and Russell~W. Schwab, \emph{Integro-differential
  equations with nonlinear directional dependence}, Indiana Univ. Math. J.
  \textbf{63} (2014), no.~5, 1467--1498.

\bibitem{KS1}
N.~V. Krylov and M.~V. Safonov, \emph{An estimate for the probability of a
  diffusion process hitting a set of positive measure}, Dokl. Akad. Nauk SSSR
  \textbf{245} (1979), no.~1, 18--20.

\bibitem{KS2}
\bysame, \emph{A property of the solutions of parabolic equations with
  measurable coefficients}, Izv. Akad. Nauk SSSR Ser. Mat. \textbf{44} (1980),
  no.~1, 161--175, 239.

\bibitem{Landkof}
N.S. Landkof, \emph{{Foundations of Modern Potential Theory}}, Graduate Studies
  in Math., Springer, 1972.

\bibitem{L}
Yan~Yan Li, \emph{Remark on some conformally invariant integral equations: the
  method of moving spheres}, J. Eur. Math. Soc. (JEMS) \textbf{6} (2004),
  no.~2, 153--180.

\bibitem{LiebLoss}
E.H. Lieb and M.~Loss, \emph{{Analysis}}, 2 ed., Graduate Studies in Math.,
  AMS, 2001.

\bibitem{LiebYau}
Elliott~H. Lieb and Horng-Tzer Yau, \emph{The stability and instability of
  relativistic matter}, Comm. Math. Phys. \textbf{118} (1988), no.~2, 177--213.
  \MR{956165}

\bibitem{MO}
Mathoverflow, \emph{{Can I approximate Schwartz functions which integrate to
  zero by $C_0^\infty$ functions which integrate to zero?}},
  http://mathoverflow.net/q/167383.

\bibitem{Moser}
J{\"u}rgen Moser, \emph{A new proof of {D}e {G}iorgi's theorem concerning the
  regularity problem for elliptic differential equations}, Comm. Pure Appl.
  Math. \textbf{13} (1960), 457--468.

\bibitem{MuehlemannThesis}
A.\ Muehlemann, \emph{{Regularity of eigenfunctions of Schr\"odinger operators
  with $L^p$\-potential}}, Master thesis, LMU Munich, 2011, url:
  www.mathematik.uni-muenchen.de/\textasciitilde lerdos/Stud/Muhleman.pdf.

\bibitem{Nash}
J.~Nash, \emph{{Parabolic Equations}}, e Naturali \textbf{3} (1957), 754--758.

\bibitem{ReedSimonII}
M.~Reed and B.~Simon, \emph{{Fourier Analysis, Self-Adjointness}}, Academic
  Press, 1975.

\bibitem{Ros-OtonSerra14}
X.~Ros-Oton and J.~Serra, \emph{{The Dirichlet problem for the fractional
  Laplacian: Regularity up to the boundary}}, J. Math. Pures Appl. \textbf{101}
  (2014), 275--302.

\bibitem{RudinFA}
W.~Rudin, \emph{{Functional Analysis}}, McGraw-Hill, 1991.

\bibitem{Silvestre07}
L.~Silvestre, \emph{{Regularity of the obstacle problem for a fractional power
  of the Laplace operator}}, Comm.\ Pure Appl.\ Math.\ \textbf{60} (2007),
  67--112.

\bibitem{SteinSIDPF}
E.M. Stein, \emph{{Singular Integrals and Differentiability Properties of
  Functions}}, Princeton Univ. Press, 1970.

\bibitem{SteinWeiss}
E.M. Stein and G.~Weiss, \emph{{Introduction to Fourier Analysis on Euclidean
  Spaces}}, Princeton Univ. Press, 1971.

\bibitem{Vladimirov}
V.~S. Vladimirov, \emph{Equations of mathematical physics}, Translated from the
  Russian by Audrey Littlewood. Edited by Alan Jeffrey. Pure and Applied
  Mathematics, vol.~3, Marcel Dekker, Inc., New York, 1971.

\bibitem{Yudovic61}
V.I. Yudovic, \emph{{Some estimates connected with integral operators and with
  solutions of elliptic equations \textnormal{(Russian)}}}, Dokl. Akad. Nauk
  SSSR \textbf{138} (1961), 805--808.

\end{thebibliography}

\end{document}